\documentstyle[12pt]{article}

\makeatletter
\renewcommand{\@biblabel}[1]{#1.}
\makeatother

\newtheorem{De}{\sc \hspace{0.4cm} Definition}[section]
\newtheorem{Ex}{\sl \hspace{0.4cm} Example}[section]
\newtheorem{Rem}{\sl \hspace{0.4cm} Remark}[section]
\newtheorem{Sa}{\sc \hspace{0.4cm} Proposition}[section]
\newtheorem{Le}[Sa]{\sc \hspace{0.4cm} Lemma}
\newtheorem{Th}[Sa]{\sc \hspace{0.4cm} Theorem}

\newcommand{\RR}{{\rm I\kern-0.14em R }}
\newcommand{\NN}{{\rm l\kern-0.14em N}}
\newcommand{\CC}{{\rm\raise 0.192ex\vbox{\hrule height
                  1.22ex width 0.8pt}\kern-0.29em C}}
\newcommand{\ZZ}{ {\sf Z}\hspace{-0.4em}{\sf Z}\ }

\begin{document}
%$\mathbb{R}$
%\[ A\  \substack{ \hookrightarrow \\ d}\ B  \]
\begin{center}
{\sf \Large
 A  Quantum Deformation of Invariants
 of  Higher Binary Forms}
\end{center}

%\noindent
%short title:

\begin{center}
Frank Leitenberger
\end{center}

\begin{flushleft}
{\large FRANK LEITENBERGER } \\
{\small \it Fachbereich Mathematik, Universit\"at Rostock,
Rostock, D-18051, Germany. \\
e-mail: frank.leitenberger@mathematik.uni-rostock.de }
\end{flushleft}

%\begin{flushleft}
%02.20.Qs; 03.65.Fd
%\end{flushleft}

%\renewcommand{\baselinestretch}{2}
\small \normalsize

\begin{abstract}
We use the theory of the quantum group $U_q(gl(2,\RR ))$
in order to develop a quantum theory of invariants
and show a decomposition of invariants into a Gordan-Capelli series.
Higher binary forms are introduced on the basis of braided algebras.
We define quantised invariants and give basic examples.
We show that the symbolic method of Clebsch and Gordan
works also in the quantised case.
We discuss the deformed discriminant
of the quadratic and the cubic form,
the deformed invariants $I_1$, $I_2$ of the quartic form
and further invariants  without a classical analog. \\
%{\it Key Words:} invariant theory; binary forms; quantum groups.
\end{abstract}

%\newpage
\vspace{0.5cm}
\setcounter{equation}{0}

\begin{center}
\Large 1. INTRODUCTION
\end{center}
\setcounter{section}{1}

In this article we develop a theory of a quantum deformation
of the invariant theory of higher binary forms.
The invariant theory was developed by Cayley, Sylvester,
Clebsch, Gordan, Capelli etc.
(cf. \cite{Cl,Go,Sa}).
The invariant theory is one of the historical roots
of the modern representation theory,
which has a well developed quantum deformation.
For the relation between both theories
we refer to \cite{Di,We}.

To find quantum deformations
of some basic notions of algebraic geometry
(e.g. curves of higher genus) is an open problem.
Plane curves with genus $g\geq 1$ can be described
by homogeneous ternary forms of degree $k\geq 3$.
These ternary forms are related to binary forms
by the "\"Ubertragungsprinzip", cf. \cite{C2}.
This is a motivation to consider
a quantum deformation of an invariant theory
of the higher binary forms.

We will demonstrate
that the basic concepts of the theory
of the quantum group  $U_q (sl(2,\RR ))$ can be used,
to build a quantum deformation of the invariant theory
of the higher binary forms.
Our approach is based on replacing
the commutative algebra of classical homogeneous
coordinates $x_i,y_i$ of $n$ points of the real line
by a braided algebra.
If the quantum parameter $q$ is not a root of unity,
the noncommutative subalgebra of invariant elements is generated
by a deformation of the bracket symbols $(ij)$ of Clebsch.
We show that there is a decomposition of invariants
into Gordan-Capelli series.

Using the braided algebra we define the notions of $n$-forms
and invariants of $n$-forms. The $n$-forms are related to
quantised polynomials. We define quantum analogs
of symmetric functions and  power sums and demonstrate Newton relations.

It turns out that the Clebsch-Gordan symbolic method
also works in our situation.
For this purpose we introduce the notions of universal forms
and universal invariants, and we consider quantum $n$-forms
and invariants as realisations of these forms.

The Clebsch-Gordan method works even in the case
when $q$ is a root of unity. But in this case
we cannot represent all
invariant elements of the algebra of noncommutative coordinates
by symbols $(ij)$.
We find that certain coefficients play the role of invariants.

By computer calculations we derive the simplest  invariants
of linear, quadratic, cubic and quartic forms.
It turns out that the ring of invariants has a richer structure
than in the classical case. For example, we obtain a
nonvanishing invariant for the linear form, an invariant of
degree three for the quadratic form and an invariant of
degree two for the cubic form,
which have no analogs in
the classical case.

Furthermore we consider
deformations of discriminants.
For the quadratic form
the situation is similar to the classical case.
For the cubic form we obtain a four-dimensional space
of invariants with properties of the discriminant.
One of them plays a special role.
For the quartic form we encounter an obstruction
for a representation of the discriminant
by the basic invariants of degree two and three.

We used the computer algebra programs Mathematica 3.0
and FELIX forming the computations.

\begin{center}
\Large 2. PRELIMINARIES
\end{center}

\setcounter{De}{0}
\setcounter{Ex}{0}
\setcounter{Rem}{0}
\setcounter{Sa}{0}

In this section we introduce some basic concepts
about the quantum group $U_q (gl (2,\RR ))$
and the braided module algebras $H_I$.

We will use the $q$-numbers
\[  [i]_{q^2} = [i] := q^{i-1} + q^{i-3} + ... + q^{-i+1}
= \frac{ q^i - q^{-i} }{ q - q^{-1} } \]
and
\[  (i)_{q^2} = (i) := 1 + q^2 + q^4 + ... + q^{2i-2}
= \frac{ q^{2i} - 1 }{ q^2 - 1 } \]
for $i\in\ZZ$, i.e. $(i) = q^{i-1} [i].$

Let $U_q (gl(2,\RR ))$, $|q|=1$,
$q\neq \pm 1$, $\pm i$
be the unital $*$-algebra, determined by generators
$E$, $F$, $K$, $L$, $K^{-1}$, $L^{-1}$,
where $L$, $L^{-1}$ are central
and relations
\[    LL^{-1}= L^{-1}L   =1,
 \ \  KK^{-1}= K^{-1}K   =1,                   \]
\[    KE     = q      EK,
 \ \  KF     = q^{-1} FK,                      \]
\[   [E,F]=\frac{ K^2-K^{-2} }{ q-q^{-1} }.    \]
We endow $U_q (gl(2,\RR ))$ with the involution
\[  E^*=F,     \ \ \ \ \ F^*=E,      \ \ \ \ \
    K^*=K^{-1},\ \ \ \ \ L^*=L^{-1}.         \]

For an arbitrary ordered index set $I$ we consider the
unital $U_q(gl(2,\RR ))$-module algebra $A_I$,
which is freely generated by the variables
$x_i,y_i$,  $i\in I$. To determine
the action of $U_q(gl(2,\RR ))$ on $A_I$, we set
\begin{eqnarray}
\begin{array}{lll}
L 1   =  1                   ,\ \ \ \ \   &
L x_i =  q^{ \frac{1}{2}} x_i,\ \ \ \ \   &
L y_i =  q^{ \frac{1}{2}} y_i,\\
K 1   =  1                   ,\ \ \ \ \   &
K x_i =  q^{-\frac{1}{2}} x_i,            &
K y_i =  q^{ \frac{1}{2}} y_i,\\
E 1   =  0                   ,\ \ \ \ \   &
E x_i =  q^{\frac{1}{2}}  y_i,            &
E y_i =  0                   ,\\
F 1   =  0                   ,\ \ \ \ \   &
F x_i =  0                   ,            &
F y_i =  q^{-\frac{1}{2}} x_i
\end{array}
\end{eqnarray}
and require that

\begin{eqnarray}
\begin{array}{lll}
 K(ab)  &  =  &  K(a) K(b),               \\
 E(ab)  &  =  &  E(a) K(b)+K^{-1}(a)E(b), \\
 F(ab)  &  =  &  F(a) K(b)+K^{-1}(a)F(b)
\end{array}
\end{eqnarray}
for $a,b \in A_I$.
A proof can be given along the lines of \cite{K}, p. 19,
where a similar module occurs.

We note that
\begin{eqnarray}
\begin{array}{llc}
E(x_i^m y_i^n) & = & q^\frac{2-m+n}{2} [m] x_i^{m-1} y_i^{n+1}, \\
F(x_i^m y_i^n) & = & q^\frac{m-n}{2} [n] x_i^{m+1} y_i^{n-1}.
\end{array}
\end{eqnarray}

Furthermore we consider the ideal $J$ of $A_I$
which is generated by the elements
\begin{eqnarray*}
  x_iy_i  -     q   y_i x_i, \ \ \ \ \
  x_jx_i  -     q^2 x_i x_j, \ \ \ \ \
  y_jy_i  -     q^2 y_i y_j,
\end{eqnarray*}
\begin{eqnarray}
  x_jy_i  -     q   y_i x_j
       - ( q^2 -1 ) x_i y_j, \ \ \ \ \
  y_jx_i  -     q   x_i y_j
\end{eqnarray}
for $i<j$.

One can show that the ideal $J$ is $U_q(gl(2,\RR ))$-invariant
(i.e. $E(J)=F(J)=0$, $K(J)\subseteq J$ and $L(J)\subseteq J$)
(cf. \cite{Ma}).

Therefore the action of $U_q(gl(2,\RR ))$ on $A_I$ induces an action on
$H_I:=A_I/J$.
In the following we identify $x_i$ with its image under the quotient map
$ A_I \rightarrow H_I $.

I.e. $H_I$ is the unital $U_q(gl(2,\RR ))$-module algebra
with generators $x_i,y_i$ $i\in I$ and relations
\begin{eqnarray}
\begin{array}{lll}
  x_iy_i             &  = &      q   y_i x_i,   \\
  x_jx_i             &  = &      q^2 x_i x_j,   \\
  y_jy_i             &  = &      q^2 y_i y_j,   \\
  x_jy_i             &  = &      q   y_i x_j
                  + ( q^2 -1 )       x_i y_j,   \\
  y_jx_i             &  = &      q   x_i y_j
\end{array}
\end{eqnarray}
for $i<j$.
Since the relations (5) are $*$-invariant for the involution
$x_i^*=x_i$, $y_i^*=y_i$ we consider $H_I$ as a $*$-module
algebra.
If $I=\{ 1,2,...,n \}$ we denote the algebra by $H_n$.

From (5) we derive the relations
\begin{eqnarray}
\begin{array}{l}
x_j y_i^k = q^k y_i^k x_j + ( q^{2k}-1 ) x_i y_i^{k-1} y_j,  \\
x_j^k y_i = q^k y_i x_j^k + ( q^{2k}-1 ) x_i x_j^{k-1} y_j.
\end{array}
\end{eqnarray}

We say that that an element $a$  of the module $H_I$ is
{\it  homogeneous of degree $k$} if $L a = q^\frac{k}{2} a$.

\begin{Rem}  \rm
$H_I$ is a {\it braided} module algebra  (i.e.
$(x_i+x_j)(y_i+y_j)=q (y_i+y_j)(x_i+x_j)$ for $i\neq j$),
(cf. \cite{Ma}). Furthermore the algebra $H_I$ has the
PBW property for an arbitrary order of the variables
$x_i,y_i$, $i\in I$.
We note that $H_I$ is up to reordering the points
the only module algebra with quadratic relations
which is braided and has the PBW property.
\end{Rem}

\begin{Rem} \rm
We can consider the Eqs. (5) as rules
in order to express an element $a\in H_I$ by the PBW basis.
These rules do not change the type of
homogenity with respect to the indices $i\in I$.
\end{Rem}

At the end of this section we extend $H_I$
to a division algebra $Q_I$. Let $R$ be an algebra,
$\sigma$ an endomorphism of $R$
and $\delta$ a $\sigma$-derivation with
\[ \delta (ab) = \delta (a) b+\sigma (a) \delta (b).\]
We say that $S=R[x;\sigma, \delta ]$
is an {\it Ore extension of $R$}
if $S$ is freely generated over $R$ by an element $x$
subject only to the relation
\[     xa = \sigma (a) x  +  \delta (a)   \]
for $a\in R$ (cf. \cite{Mc}).

\begin{Le}
The algebra $H_I$ is an iterated Ore extension.
\end{Le}

{\it Proof.}
We restrict the consideration to the finite algebras $H_n$.
In the case of an infinite index set $I$ one can proceed by
transfinite induction.
We consider the tower of subalgebras
\[  H_0=\CC \subset H_1' \subset H_1 \subset \ \cdots \
            \subset H_n' \subset H_n   \]
where $H_i'$ is the subalgebra generated by $H_{i-1}$ and $y_i$.

We have $H_i' \cong H_{i-1} [y_{n+1},\sigma,\delta]$
where $\delta = 0$ and $\sigma$ is determined by
\[  \sigma (x_i)=q   x_i,\ \ \ \ \
    \sigma (y_i)=q^2 y_i. \]
Furthermore we have $H_i \cong H_i' [x_{n+1},\sigma,\delta]$
where $\sigma$ and $\delta$ are determined by
\[   \sigma (x_i)    = q^2   x_i,\ \ \ \ \
     \sigma (y_i)    = q     y_i,\ \ \ \ \
     \sigma (y_{n+1})= q     y_{n+1},      \]
\[   \delta (x_i)    = 0,\ \ \ \ \
     \delta (y_i)    = (q^2-1) x_i y_{n+1},\ \ \ \ \
     \delta (y_{n+1})= 0. \]
The verification of the isomorphies
is similar to \cite{Ka}, p. 81.
The assertion follows.
${\bullet}$

\begin{Sa}
The algebra $H_I$ has no zero divisors
and $H_I$
has an uniquely determined right quotient division ring $Q_I$,
whose elements are the right quotients
$ab^{-1}$ for $a,b\in R$, $b\neq 0$.
\end{Sa}

{\it Proof.}
Ore proved that a noncommutative ring without zero divisors
has an uniquely determined quotient division ring if
\begin{eqnarray}
  aR \cap bR \neq \{ 0\}\ \ \ \ \ \ \ \forall a,b \in R,
                            \ \ \ \ \ a,b\neq 0,
\end{eqnarray}
i.e., two nonvanishing elements $a,b$
have a common right multiple (cf. \cite{O}).
Curtis proved that an iterated Ore extension
of a ring without zero divisors with the property (7)
is again a ring without zero divisors with the property (7)
(cf. \cite{Cu}).
The assertion follows from both results together with Lemma 2.1.
${\bullet}$

\begin{Ex} \rm
For $x_2$ and $y_1$ we find the common right multiple
\[    x_2 ( x_2 y_1 - (q^2-q^{-2}) x_1 y_2 )
    = y_1 ( q^2 x_2^2 ).  \]
Therefore we can represent the element $y_1^{-1} x_2$
by the right quotient
$ q^2 x_2^2 ( x_2 y_1 - (q^2-q^{-2}) x_1 y_2 )^{-1}$.
Note, that in general it is not possible
to find a common right multiple of $a$ and $b$ of degree
$\leq deg(a)+deg(b)$.
\end{Ex}

\begin{Rem} \rm
We can interprete the variables $x_i, y_i$
as noncommutative homogeneous coordinates
of different points of the projective real line (cf. \cite{Le}).
%We adjoin elements $y_i^{-1}$, $i\in H_I$ with
%$y_i y_i^{-1} = y_i^{-1} y_i = 1$ to $H_I$.
The noncommutative real coordinate is given by
\[  v_i: = q^{\frac{1}{2} } x_i y_i^{-1} \in Q_I. \]
We have $v_i^*=v_i$ and
\[ v_j v_i  = q^2 v_i v_j + (1-q^2) v_i^2
\ \ \ \ \ {\rm for} \ \ \ \ \   i<j. \]
For $q=-1$ we obtain a noncommutative algebra of homogeneous coordinates
for the classical commutative algebra of real coordinates.
\end{Rem}

\begin{center}
\Large 3. THE INVARIANTS OF THE ALGEBRA $H_I$
\end{center}
\setcounter{section}{3}
\setcounter{De}{0}
\setcounter{Ex}{0}
\setcounter{Rem}{0}
\setcounter{Sa}{0}

\begin{De} \rm
We say that the element $a$ of $H_I$ is an
{\it invariant element of degree $k$},
if $Ea=Fa=0$, $Ka=a$ and $La=q^\frac{k}{2}a$.
\end{De}

By the rules (2) sums and products of invariant elements
are again invariant elements.
Therefore they form a subalgebra
$H_I^{Inv} \subset H_I$.

\begin{Ex} \rm
The simplest invariant element of $H_I$ is given by
\begin{eqnarray*}
(ij):= q^{-\frac{1}{2}} x_i y_j- q^{\frac{1}{2}} y_i x_j,
\ \ \ \ \    i,j \in I.
\end{eqnarray*}
We have
\[ (ji)  = - (ij),\ \ \forall i,j, \]
\[ (ij)^*=   (ij),\ \ \forall i,j  \]
and
\[ (ii)  =    0,  \ \ \forall i.   \]
The condition $(ji)  = -(ij)$
is equivalent to
the braiding condition \\
$ (x_i+x_j) (y_i+y_j) = q (y_i+y_j) (x_i+x_j) $
for $i\neq j$, cf. Remark 2.1.
\end{Ex}

\begin{Le} Let $i<j<k$  be three ordered indices. Then the
following commutation relations hold:
\begin{eqnarray}
\begin{array}{lllllllll}
(jk) x_i  & = & q^3 &  x_i (jk),  & \ \ \ \ \ &
(jk) y_i  & = & q^3 &  y_i (jk),  \\
(jk) x_j  & = & q   &  x_j (jk),  & \ \ \ \ \ &
(jk) y_j  & = & q   &  y_j (jk),  \\
x_j (ij)  & = & q   &  (ij) x_j,  & \ \ \ \ \ &
y_j (ij)  & = & q   &  (ij) y_j,  \\
x_k (ij)  & = & q^3 &  (ij) x_k,  & \ \ \ \ \ &
y_k (ij)  & = & q^3 &  (ij) y_k.
\end{array}
\end{eqnarray}
\end{Le}

\begin{Sa}
(i) Let $i<j<k<l$  be four ordered indices. Then
the following commutation relations hold:
\begin{eqnarray}
\begin{array}{lll}
 (kl)(ij) & = & q^{ 6} (ij)(kl),   \\
 (jk)(il) & = &        (il)(jk),   \\
 (jl)(ik) & = & q^{ 4} (ik)(jl) + (q^{ 4}-q^{ 6}) (ij)(kl), \\
 (jk)(ij) & = & q^{ 4} (ij)(jk),   \\
 (ik)(ij) & = & q^{ 2} (ij)(ik),   \\
 (jk)(ik) & = & q^{ 2} (ik)(jk).
\end{array}
\end{eqnarray}

\noindent
(ii) Let $i,j,k,l$ be four different indices. Then
the three invariant elements
\[   (ij)(kl),\ \  (ik)(jl),\ \  (il)(kj) \]
commute.

\noindent
(iii) Let $i<j<k<l$  be four ordered indices.
Then we have the identities
\begin{eqnarray}
 q^4 (ij)(kl) + q^2 (ik)(lj) + (il)(jk)=0
\end{eqnarray}
(Grassmann-Pl\"ucker relations).
\end{Sa}

A way to check the identities of Lemma 3.1 and Proposition 2.2
is indicated in \cite{Le}, p. 816.

\begin{Rem} \rm
We can rewrite Proposition 3.2 (ii), (iii) in a more
symmetric formulation. Let
\[ Z_1:=q^3 (ij)(kl),                 \]
\[ Z_2:=q^2 (ik)(lj)+(q^4-q^3)(ij)(kl), \]
\[ Z_3:=    (il)(jk).                   \]
According to Proposition 3.2 the elements
$Z_1,\ Z_2,\ Z_3$ commute, are $*$-invariant and we have
\[ Z_1+Z_2+Z_3 = 0. \]
\end{Rem}

\begin{Rem} \rm
If we interprete $x_i, y_i$ as noncommutative
projective coordinates, then the six elements
$-Z_1 Z_2^{-1},\ -Z_1 Z_3^{-1},\ -Z_2 Z_1^{-1},\
 -Z_2 Z_3^{-1},\ -Z_3 Z_1^{-1},\ -Z_3 Z_2^{-1}\ \in Q_I\ $
get the meaning of the six possible cross ratios of the four
zeros $v_i:=q^\frac{1}{2} x_i y_i^{-1}$ of the form $f$.
For example, we have
\[     - Z_1 Z_3^{-1} =
   \frac{1}{q} (v_i-v_j)(v_j-v_k)^{-1} (v_k-v_l)(v_i-v_l)^{-1} \]
(cf. \cite{Le}).
\end{Rem}

\begin{Le}
Let $i<j<k$.
The following formulas hold:

\noindent
(i)
\[ (ij)(ik)      =      x_i^2        y_j y_k
                   -    x_i y_i (q   x_j y_k + y_j x_k )
                   +q^2 y_i^2        x_j x_k,          \]

\[ (ik)(jk)      =          x_i x_j           y_k^2
                   - (q x_i y_j +   y_i x_j ) x_k y_k
                   +q^2     y_i y_j           x_k^2,   \]

\[ (ij)(jk)      =  \frac{1}{q^2} x_i x_j y_j y_k
                   - ( x_i y_j y_j x_k +   y_i x_j x_j y_k)
                   +         q    y_i x_j y_j x_k,   \]

\[ (ij)(ij)      =       x_i^2 y_j^2
                   - [2] x_i y_i x_j y_j
                   + q^2 y_i^2 x_j^2   \]
\[             =            \frac{ 1 }{q^2}  x_j^2 y_i^2
                       -    \frac{[2]}{q^2}  x_j y_j x_i y_i
                       +                     y_j^2   x_i^2, \]

\noindent
(ii)
\[ (ij)^n = \sum_{k=0}^{n}
        \left[  \begin{array}{c} n \\ k \end{array}   \right]
    (-1)^k     q^{\frac{1}{2} n^2 +(k-n)(1+k)}
              x_i^{n-k} y_i^k x_j^{k} y_j^{n-k}, \]

\noindent
(iii)
\begin{eqnarray*}
\begin{array}{ll}
 (ij)(ik)(jk)=  &
    q^{-\frac{1}{2}}  x_i x_i x_j y_j y_k y_k    -
    q^{ \frac{3}{2}}  x_i x_i y_j y_j x_k y_k    -
    q^{ \frac{3}{2}}  x_i y_i x_j x_j y_k y_k     \\
 &
  - q^{-\frac{3}{2}}  x_i y_i x_j y_j x_k y_k    +
    q^{ \frac{5}{2}}  x_i y_i x_j y_j x_k y_k    +
    q^{ \frac{3}{2}}  x_i y_i y_j y_j x_k x_k     \\
 &
  + q^{ \frac{3}{2}}  y_i y_i x_j x_j x_k y_k    -
    q^{ \frac{7}{2}}  y_i y_i x_j y_j x_k x_k.
\end{array}
\end{eqnarray*}

\noindent
(iv) Let $u_i^{(0)}:=x_i$ and $u_i^{(1)}:=y_i$. We have
\[
(1,2n)(2,2n-1)(3,2n-2) \cdots (n,n+1)
\]
\[
=
\sum_{ i_1,i_2, \cdots , i_n = 0 }^1
q^{ \frac{3}{2} n(n-1) }  (-q)^{ i_1+i_2+ \cdots +i_n }
u_{1  }^{(i_1)}
u_{2  }^{(i_2)}
\cdots
u_{n  }^{(i_n)}
u_{n+1}^{(1-i_n)}
u_{n+2}^{(1-i_{n-1})}
\cdots
u_{2n }^{(1-i_1)}.
\]
\end{Le}
Here \[ \left[  \begin{array}{c}
          n \\ i \end{array}   \right]:=
        \frac{ [n]! }{ [n-i]! [i]! },   \]
denotes the $q$-binomial coefficient, where
\[   [n]! := [1] [2] ... [n],      \]
$i,n\in \NN $.
(If $q^{2d} = 1$ for all $d\in \NN$, but $q^2 \neq 1$,
we define
$\left[  \begin{array}{c}
          n \\ i \end{array}   \right]$
by a limit.)

A proof of Lemma 3.3 can be given by
an explicite calculation or by induction.

\begin{flushleft}
3.1 \it  Polarisation operators
\end{flushleft}

We introduce the polarisation operators
$ \Delta_{kl}: H_I  \rightarrow  H_I $, $k\neq l$.
We fix a PBW basis, which consists of the elements
\[   x_{a_1}^{ i_1 }   y_{a_1}^{ j_1 }
     x_{a_2}^{ i_2 }   y_{a_2}^{ j_2 }  ...
     x_{a_k}^{ i_k }   y_{a_k}^{ j_k }       \]
with $a_1 < a_2 < ... < a_k$.
We set
\begin{eqnarray*}
\Delta_{kl} ( x_k ) = x_l,\ \ \ \ \ \ \ &
\Delta_{kl} ( y_k ) = y_l
\end{eqnarray*}
and
\begin{eqnarray*}
\Delta_{kl} ( x_i ) = 0,  \ \ \ \ \ &
\Delta_{kl} ( y_i ) = 0
\end{eqnarray*}
for   $i\neq k$ and require
\[ \Delta_{kl} ( ab ) =    \Delta_{kl} (a)\ b +
                        a\ \Delta_{kl} (b),     \]
if the expressions $a$, $b$ and $ab$
are elements of the PBW basis.

If $f$ is a polynomial,
which is homogeneous with respect to $k\in I$ of degree $n_k$
we have the {\it Euler identity}
\[     \Delta_{kk} f = n_k\ f.  \]

\begin{Rem} \rm
One obtains different linear operators $\Delta_{kl}$
for different PBW bases. For our purposes it is sufficient
to work with the above fixed basis.
\end{Rem}

Furthermore, if $q^{2d}\neq 1$ for all $d\in \NN$,
we define a linear operator $P_{kl}$ by
\[ P_{kl}  a := \frac{1}{(n_k)_{q^2}  } \Delta_{kl} a, \]
if $a$ is a monomial of $H_I$,
which is homogeneous
with respect to the index $k$ of degree $n_k$.

\begin{Sa} Let $q^{2d} \neq 1$ for all $d\in \NN$.

\noindent
(i) We have $P_{ij} (ik)=(jk)$, i.e. $P_{ij} (ij) =0$.

\noindent
(ii)  For $k<l$
and $ \beta \in H_{ \{ i\in I|i>l \} }$
we have
\[  P_{kl} \left(
    x_k^{m-i} y_k^{i} y_l^{n-j} x_l^{j} \beta \right) =  \]
\[  =   \frac{1}{[m]} (
            [m-i] x_k^{m-i-1} y_k^{i}   x_l^{n-j+1} y_l^{j} \beta
+ q^{m-i-n+j} [i]
   x_k^{m-i} y_k^{i-1} x_l^{n-j}   y_l^{j+1} \beta ) \]
and
\[  P_{lk} \left(
    x_k^{m-i} y_k^{i} x_l^{n-j} y_l^{j} \beta \right) = \]
\[  =   \frac{1}{[n]} (
    q^{j-i} [n-j]  x_k^{m-i+1} y_k^{i}   x_l^{n-j-1} y_l^{j}\beta
  +         [j]
   x_k^{m-i}   y_k^{i+1} x_l^{n-j}   y_l^{j-1}\beta ) .\]

\noindent
(iii) $P_{kl}$ commutes with the action of $U_q(gl(2,\RR )$,
i.e., if $ f \in  H^{Inv}_I $,
      then $P_{kl} (f) \in H^{Inv}_I $.

\end{Sa}

{\it Proof.}
(i) For  $i<k$ we have
\[ P_{ij} (ik) = \Delta_{ij} q^{-\frac{1}{2}}(x_i y_k-q y_i x_k)
  =  q^{-\frac{1}{2}} (  x_j y_k - q y_j x_k ) = (jk).        \]
The case $i>k$ is analogous.

\noindent
(ii)  From the definition of $\Delta_{kl}$ we obtain
      $ \Delta_{kl} (x_k^m)= (m) x_k^{m-1} x_l $ and
      $ \Delta_{kl} (y_k^m)= (m) y_k^{m-1} y_l $.
It follows
\[  \Delta_{kl} ( x_k^{m-i}   y_k^{i}   x_l^{n-j}   y_l^{j} )
 =  (m-i)  x_k^{m-i-1} x_l   y_k^{i}  x_l^{n-j} y_l^{j} +
    (  i)  x_k^{m-i}   y_k^{i-1} y_l  x_l^{n-j} y_l^{j}.
\]
Using (5) and (6) we find
\[ =(m-i) q^i         x_k^{m-i-1} y_k^{i}   x_l^{n-j+1} y_l^{j} +
    (  i) q^{2m-2i-n+j} x_k^{m-i}   y_k^{i-1} x_l^{n-j}   y_l^{j+1}.\]
Since $ (m) = q^{m-1} [m] $ we have
\[ =[m-i] q^{m-1}         x_k^{m-i-1} y_k^{i}   x_l^{n-j+1} y_l^{j} +
    [  i] q^{2m-i-n+j-1}  x_k^{m-i}   y_k^{i-1} x_l^{n-j}   y_l^{j+1}. \]
The formula for $P_{kl}$ follows.
The proof for $P_{lk}$ is analogous.

\noindent
(iii) By the definition of $K, L$ and $\Delta_{kl}$ we have
$K \Delta_{kl} f = \Delta_{kl} K f$ and
$L \Delta_{kl} f = \Delta_{kl} L f$.

From (ii) and (3) it follows that
\[ \Delta_{kl} E ( x_k^m y_k^n )
  = q^\frac{2-m+n}{2} [m] \Delta_{kl}  ( x_k^{m-1} y_k^{n+1} )        \]
\[= q^\frac{2-m+n}{2} [m]
([m-1] x_k^{m-2} y_k^{n+1} x_l + q^{m-1} [n+1] x_k^{m-1} y_k^{n} y_l )
\]
and
\[ E \Delta_{kl} ( x_k^m y_k^n ) =
       [m] E      (x_k^{m-1} y_k^n    ) K(x_l)
 +     [m] K^{-1} (x_k^{m-1} y_k^{n})   E(x_l)
 + q^m [n] E      (x_k^{m}   y_k^{n-1}) K(y_l)  \]

\[ = q^\frac{2-m+n}{2}   [m][m-1]  x_k^{m-2}  y_k^{n+1}  x_l
   + q^\frac{m-n  }{2}   [m]       x_k^{m-1}  y_k^n      y_l
   + q^\frac{2+m+n}{2}   [n][m]    x_k^{m-1}  y_k^n      y_l  \]
\[  =  q^\frac{2-m+n}{2}   [m]
([m-1] x_k^{m-2} y_k^{n+1} x_l+q^{m-1} [n+1] x_k^{m-1} y_k^{n} y_l)
\]
(since $q^{-n} [1]+q [n]=[n+1]$).
Therefore we have
$E \Delta_{kl} f_k = \Delta_{kl} E f_k$
for $f_k=x_k^m y_k^n$.

Now let $f$ be a PBW basis element.
We consider the decomposition $f=f_- f_k  f_+$ with
PBW basis elements $f_+ \in H_{ \{ i\in I|i<k \} }$,
                   $f_- \in H_{ \{ i\in I|i>k \} }$. We have
\[   E \Delta_{kl}   f             =
     E \Delta_{kl}  ( f_- f_k  f_+ ) =
     E ( f_- \Delta_{kl} (f_k) f_+ )          \]
\[   =
   E      (f_-) (\Delta_{kl} (K      (f_k))) K(f_+)+
   K^{-1} (f_-) (\Delta_{kl} (E      (f_k))) K(f_+)+
   K^{-1} (f_-) (\Delta_{kl} (K^{-1} (f_k))) E(f_+) \]
\[   =  \Delta_{kl} E ( f_- f_k  f_+ )
     =  \Delta_{kl} E  f.                          \]
Similarly we find  $F \Delta_{kl} f = \Delta_{kl} F f$
for all elements $f$ of the PBW basis.
I.e., $E,F,K$ commute with $\Delta_{kl}$ and therefore
with $P_{kl}$. The proposition
follows.$\bullet$

\begin{Sa}
Let $q^{2d} \neq 1$ for all $d\in \NN$
and let $f\in H_I$ be a linear combination of monomials
\[  x_k^{m-i} y_k^{i} x_l^{n-j} y_l^{j} \beta_{ij} \]
with $k<l$, $m,n\geq 1$
and $\beta_{ij} \in H_{ \{ i\in I|i>l \} }$.

\noindent
(i) Then we can represent $f$ in the unique form
\[   f =   P_{lk} P_{kl} f + (kl) a \]
or
\[   f =   P_{kl} P_{lk} f + (kl) b \]
with $a,b \in H_I$ and  $b=\frac{ [m][n+1] }{ [m+1][n] } a$.

\noindent
(ii) We have
\[ [m][n+1] P_{lk} P_{kl} f -
   [n][m+1] P_{kl} P_{lk} f = [m-n] f  \]

\noindent
(iii) If  $f\in H_I^{Inv}$, then $a\in H_I^{Inv}$.
\end{Sa}

{\it Proof.}
(i) It is sufficient to consider the case,
where $f$ is the monomial
$ x_k^{m-i} y_k^{i} x_l^{n-j} y_l^{j}$.
With Proposition 3.4 (ii) we find
\[  P_{lk} P_{kl} f  =  \frac{1}{[m]}
    P_{lk} (  [m-i] x_k^{m-i-1} y_k^{i}   x_l^{n-j+1} y_l^{j}
+ q^{m-i-n+j} [i]   x_k^{m-i}   y_k^{i-1} x_l^{n-j}   y_l^{j+1} )
\]
\[ =\frac{1}{[m][n+1]} (
 [m-i][n-j+1] q^{j-i} x_k^{m-i}   y_k^{i}   x_l^{n-j}   y_l^{j}
 +[m-i][j]            x_k^{m-i-1} y_k^{i+1} x_l^{n-j+1} y_l^{j-1}
\]
\[
+ [i][n-j] q^{m-n-i+j+2} x_k^{m-i+1} y_k^{i-1} x_l^{n-j-1} y_l^{j+1}
+ [i][j+1] q^{m-n-i+j}   x_k^{m-i}   y_k^{i}   x_l^{n-j}   y_l^{j}).
\]
For $m,n \geq 2$ it follows that
\[          P_{lk} P_{kl} f
         =  \frac{1}{[m][n+1]}
          x_k^{m-i-1}  y_k^{i-1}  (
   q^{m+n+2} [i]  [n-j]       x_k^2   y_l^2   + \]
\[(q^{n-2}   [m-i][n-j+1] +
   q^{m-2}   [i]  [j+1]     ) x_k y_k x_l y_l +
             [m-i][j]         y_k^2   x_l^2 )
             x_l^{n-j-1}  y_l^{j-1} \]

\[ =  \frac{1}{[m][n+1]}
          x_k^{m-i-1}  y_k^{i-1} (
       q^{ m+n-2+\frac{1}{2} } [i][n-j] x_k (kl) y_l
    -  q^{      -\frac{1}{2} } [m-i][j] y_k (kl) x_l + \]
\[   ( q^{ m+n-1 }             [i][n-j] +
       q^{    -3 }             [m-i][j] +
       q^{n-2}                 [m-i][n-j+1]+
       q^{m-2}                 [i][j+1] ) x_k y_k x_l y_l)
         x_l^{n-j-1}  y_l^{j-1}. \]
Because of  $ [i+j] =  q^{ j} [i] + q^{-i} [j]
                   =  q^{-j} [i] + q^{ i} [j] $
we obtain
\[ =  \frac{1}{[m][n+1]}
          x_k^{m-i-1}  y_k^{i-1} (
       q^{ m+n-2+\frac{1}{2} } [i][n-j] x_k (kl) y_l
    -  q^{      -\frac{1}{2} } [m-i][j] y_k (kl) x_l + \]
\[    q^{n-j+i-2} [m][n+1]   x_1 y_1 x_2 y_2)
         x_l^{n-j-1}  y_l^{j-1}.                        \]
Applying Lemma 3.1 we find
{\small
\[ = \frac{1}{[m][n+1]} (kl)
 (q^{ j-i+\frac{3}{2} } [i][n-j] x_k^{m-i} y_k^{i-1} x_l^{n-j-1} y_l^{j} -
  q^{  -m+\frac{1}{2} } [m-i][j] x_k^{m-i-1} y_k^i   x_l^{n-j} y_l^{j-1} )
  + f.  \] }
We encounter the same formula in the cases $m<2$ and $n<2$.
The first formula follows. Similarly we find
{\small
\[  P_{lk} P_{kl} f  =
  \frac{1}{[n][m+1]} (kl)
 (q^{ j-i+\frac{3}{2} }
        [i][n-j] x_k^{m-i} y_k^{i-1} x_l^{n-j-1} y_l^{j} -
  q^{  -m+\frac{1}{2} }
        [m-i][j] x_k^{m-i-1} y_k^i   x_l^{n-j} y_l^{j-1} )
  + f. \]}
The second formula and the relation
$b=\frac{ [m][n+1] }{ [m+1][n] } a$ follow.
The representation is unique,
since $H_I$ has no zero divisors (c.f. Proposition 2.2).

\noindent
(ii) The assertion is a consequence of (i)
and of the identity \\
$[m][n+1] - [n][m+1] = [m-n]$.

\noindent
(iii) By Proposition 3.4 (iii) we have
\[ 0 = Ef = P_{lk}  P_{kl} E(f) +
            E((kl)) K(a)        + K^{-1}((kl)) E(a)
          = (kl)  E(a).                 \]
Since $H_I$ has no zero divisors,
we have
$E(a)=0$. Similarly we derive $K(a)=a$ and $F(a)=0$.
$\bullet$

\begin{flushleft}
3.2 \it The main Theorem
\end{flushleft}

\begin{Th}
Let $q^{2d} \neq 1$ for all $d\in \NN $.
The subalgebra of invariant elements $H_I^{Inv}$
is generated by the elements $(ij)$, $i,j\in I $.
\end{Th}

{\it Proof.}
Since the $U_q(gl(2,\RR ))$-action on $H_I$ does
not change the type of homogenity with respect to
the indices  $i\in I$, it is sufficient to consider
invariants which are homogeneous with respect to
all indices.

Let $f \in H^{Inv}_I $ be an invariant element
of degree $n$ with $k$ series of variables
$x_{a_j},y_{a_j}$, $j=1,...,k$, $a_1<a_2< ... <a_k$,
i.e.
\[  f  =  \sum_{i_1,...,i_k}   c_{i_1,i_2,...,i_k}
            x_{a_1}^{ n_1 - i_1 } y_{a_1}^{ i_1 }
            x_{a_2}^{ n_2 - i_2 } y_{a_2}^{ i_2 } ...
            x_{a_k}^{ n_k - i_k } y_{a_k}^{ i_k }      \]
with $n_j\geq 1$, $i_j\geq 0$ and $n_1+...+n_k=n$.
We suppose $k\geq 2$ (i.e. $n\geq 2$).
For simplicity we set $a_1=1$ and $a_2=2$.

Using Proposition 3.4 (iii) we can form the $n_1$ invariants
\[                  f,\ \ \ \ \
    P_{12}          f,\ \ \ \ \
    P_{12}^2        f,\ \ \ \ \
                \cdots ,\ \ \ \ \
    P_{12}^{n_1 -1} f. \]
By Proposition 3.5 (i) we have
\begin{eqnarray*}
             f & = &  P_{21}   P_{12}   f + (12) a_0, \\
      P_{12} f & = &  P_{21}   P_{12}^2 f + (12) a_1, \\
               &   &     .....,                                  \\
      P_{12}^{n_1 -1}  f & = &
                      P_{21} P_{12}^{n_1} f + (12) a_{n_1 -1}
\end{eqnarray*}
with $a_i \in H_I^{Inv}$. It follows that
\begin{eqnarray*}
  f & = &  P_{21}   P_{12} f + (12) a_0 \\
    & = &  P_{21} ( P_{21} P_{12}^2 f
              + (12) a_1 )
              + (12) a_0                                  \\
    & = &  P_{21}^2 P_{12}^2 f
              + (12) (a_0 + P_{21} a_1) \\
                  &   &     .....                            \\
    & = &  P_{21}^{n_1} P_{12}^{n_1} f
+ (12) (a_0 +           P_{21}   a_1
               +           P_{21}^2 a_2 + ...
               +           P_{21}^{n_1 -1} a_{n_1 -1} ).
\end{eqnarray*}
Therefore we can reduce the problem for the invariant element $f$
of degree $n$ with $k$ variables
to the problems for the invariant element $P_{12}^{n_1} f$
of degree $n$ with $k-1$ variables
and the invariant element
\[ a_0 +   P_{21}   a_1
              +   P_{21}^2 a_2 + ...
              +   P_{21}^{n_1 -1} a_{n_1 -1}     \]
of degree $n-2$ with $k$ variables.
Consequently the problem reduces  to the cases $n=0,1$ or $k=1$
by a finite number of steps.

For $n=0$ we have only the invariant element $1$.

For $n=1$ there are no invariant elements,
since we infer from
\[ K( \sum_i ( \alpha_i x_i + \beta_i y_i )) =
      \sum_i (q^{-\frac{1}{2}} \alpha_i x_i +
              q^\frac{ 1}{2}   \beta_i  y_i )
    = \sum_i ( \alpha_i x_i +\beta_i y_i)     \]
that $\alpha_i=\beta_i=0$.

Finally we consider the case $k=1$.
The only elements $a\in H_I$ with $k=1$ and $K(a)=a$ are
$ a=\sum_{i=0}^n
\alpha_i x_i^m y_i^m $.
From the condition
\[  E(\sum_i \alpha_i  x_i^m y_i^m )
  =   \sum_i \alpha_i q^\frac{2-m}{2}
            [m] x_i^{m-1} y_i^{m+1} = 0 \]
it follows that $\alpha_i=0$ for $i\geq 1$,
since we supposed $q^{2d} \neq 1$,
i.e. $a=\alpha_0 1$.
Therefore we have no invariant elements for $k=1$.
The Theorem follows.
$\bullet$

\begin{flushleft}
3.3 \it  The Gordan-Capelli series
\end{flushleft}

In this section
we consider $H_I$ with $I= \{1,2, \cdots ,p   \}$.
We use Proposition 3.5 (i), (ii) in order to introduce
the quantum $\Omega$-process for the two lowest indices.

\begin{De}
Let  $ f= P_{21} P_{12} f + (12) a $.
Then we define $\Omega_{12} (f):= \frac{[n+1]}{[n]} a$.
\end{De}

By Proposition 3.5 (i) $a$ is unique and we have
\[ \Omega_{12}
   ( x_1^{m-i} y_1^{i} x_2^{n-j} y_2^{j} \beta ) \]
\[ =\frac{1}{[m][n]}
(q^{j-i+\frac{3}{2} }
     [i][n-j] x_1^{m-i} y_1^{i-1} x_2^{n-j-1} y_2^{j}\beta -
 q^{ -m+\frac{1}{2} }
 [m-i][j] x_1^{m-i-1} y_1^i  x_2^{n-j} y_2^{j-1}\beta )  \]
for $m,n \geq 1$ and $\beta \in H_{ \{i\in I| i>2 \} }$
(cf. the proof of Proposition 3.5 (i)).

Therefore we have
\begin{eqnarray}
f= P_{21} P_{12} f+\frac{ [n] }{ [n+1] } (12) \Omega_{12} (f)
\end{eqnarray}
and
\begin{eqnarray}
f= P_{12} P_{21} f+\frac{ [m] }{ [m+1] } (12) \Omega_{12} (f).
\end{eqnarray}

\begin{Rem} \rm
For $q=1$ we have
\begin{eqnarray*}
\Omega_{12} =
\frac{1}{m n}
\left|
\begin{array}{cc}
\partial_{x_1} & \partial_{y_1}    \\
\partial_{x_2} & \partial_{y_2}
\end{array}
\right|
\end{eqnarray*}
and Eq. (11)
takes the form of the {\it Capelli identity}:
\begin{eqnarray*}
\left|
\begin{array}{ll}
\Delta_{11} & \Delta_{12}    \\
\Delta_{21} & \Delta_{22}+1
\end{array}
\right| (f) =
(12)
\left|
\begin{array}{cc}
\partial_{x_1} & \partial_{y_1}    \\
\partial_{x_2} & \partial_{y_2}
\end{array}
\right|
(f)
\end{eqnarray*}
(cf. \cite{Ca,Fu}).
We have not introduced quantum deformations
of the partial derivatives
$\partial_{x_i}$, $\partial_{y_i}$.
Therefore we use Eq. (11) not as an identity,
but only for the definition of the operator $\Omega_{12}$.
This will be sufficient for the proof of Theorem 3.8.
For further applications of the Capelli identity
we refer to \cite{We,Ho}.
A Capelli identity on quantum groups
was considered in \cite{Hi} for $GL_2(\RR)$
and in \cite{No} for $GL_n(\RR)$.
\end{Rem}

The operators $P_{12}, \Omega_{12}$
have the following commutation relations

\begin{Sa}
(i)  We have
\[ P_{12} ((12)f)=q \frac{[m]}{[m+1]} (12) P_{12}f,\]
\[ P_{21} ((12)f)=q \frac{[n]}{[n+1]} (12) P_{21}f \]
for $f\in H_I$ and

\noindent
(ii)
\[ \Omega_{12} P_{12} =  \frac{1}{q}
\frac{   [n] }{ [n+1] }  P_{12} \Omega_{12}, \]
\[ \Omega_{12} P_{21} =
\frac{ q [m] }{ [m+1] }  P_{21} \Omega_{12}. \]
\end{Sa}

{\it Proof.}
(i)  A proof can be given by an explicite calculation
using Proposition 3.4 (ii).

%Let $ a = x_k^{m-i} y_k^{i} x_l^{n-j} y_l^{j} $.
%We have
%\[ (12)a= q^{m-i+j-n-\frac{1}{2}}
%       x_k^{m-i+1} y_k^{i} x_l^{n-j} y_l^{j+1} -
%       q^{m+\frac{1}{2}}
%       x_k^{m-i} y_k^{i+1} x_l^{n-j+1} y_l^{j} \]
%Using Proposition we obtain
%\[ P_{12} (12) a = \frac{1}{[m+1]}
%   ( q
%=  \]
%and
%\[ P_{12} (12) a = \frac{1}{[m]} =  \]
%The first equation follows.

\noindent
(ii) By (11) and Proposition 3.5 (ii) we have
\[                      P_{12} ((12) \Omega_{12} (f))
  =   \frac{[m]}{[m+1]} P_{12} ( f - P_{12} P_{21} f ) \]
\[=   \frac{[m]}{[m+1]}         \left( P_{12}f-
      \frac{[n][m+1]}{[m][n+1]} P_{12} P_{21} P_{12}f -
      \frac{[m-n]}{[m][n+1]}          P_{12}f \right)  \]
\[=   \frac{[n+1][m]-[m-n]}{[m+1][n+1]} P_{12}f  -
      \frac{[n]}{[n+1]}   P_{12} P_{21} P_{12}f        \]
\[=   \frac{[n]}{[n+1]} ( 1-P_{12}P_{21}) P_{12} f
=   \frac{[m][n]}{[m+1][n+1]}(12)\Omega_{12} P_{12} (f). \]
Using (i), we find
\[  q (12) P_{12} \Omega_{12} (f)) =
      \frac{[n]}{[n+1]} (12)\Omega_{12} P_{12} (f). \]
The first equation follows,
since $H_I$ has no zero divisors.
The proof of the second equation is analogous.$\bullet$

Now we demonstrate that we can represent
every invariant with $n$ variables
by powers of symbols $(12)^k$ and polars
of invariants with $n-1$ variables.

\begin{Th} Let $q^{2d}\neq 1$ for all $d\in \NN$
and let $f$ be an invariant element of degree $n$
with $k$ series of variables $x_i,y_i$, $i=1,2,...,p$,
homogeneous of degree $n_i$ with respect to $i$. Then we have
\[ f = \sum_{k=0}^{n_1}
   \frac{ \left[  \begin{array}{c} n_1 \\ k \end{array}  \right]
   \left[  \begin{array}{c} n_2 \\ k \end{array}  \right] }{
   \left[  \begin{array}{c} n_1+n_2-k+1 \\ k \end{array} \right] }
   (12)^k  P_{21}^{n_1-k} P_{12}^{n_1-k} \Omega_{12}^k  f. \]
\end{Th}

{\it Proof.} For $0\leq l \leq n_1$
we can represent $f$ in the form
\begin{eqnarray}
   f = \sum_{k=0}^{l} \alpha_k^{(l)}
   (12)^k  P_{21}^{l-k} P_{12}^{l-k} \Omega_{12}^k  f
\end{eqnarray}
with
\begin{eqnarray}
   \alpha_{k  }^{(l+1)} =
   \alpha_{k  }^{(l  )} +
   \alpha_{k-1}^{(l  )}
   \frac{ (n_2-k+1)^2 }{(n_2+l-2k+2)(n_2+l-2k+3)}
\end{eqnarray}
($\alpha_{-1}:=0$).
Since $\alpha_{0}^{(0)}=\alpha_{0}^{(1)}=1$,
the statement is true for $l=0$.
We suppose  the statement for $l$, i.e. we have
\[  f =    \sum_{k=0}^{l} \alpha_k^{(l)}
   (12)^k  P_{21}^{l-k} P_{12}^{l-k} \Omega_{12}^k  f.   \]
Applying (13) to the elements
\[  P_{12}^{l-k} \Omega_{12}^k f           \]
we find
\[  f  =    \sum_{k=0}^{l} \alpha_k^{(l)}
(12)^k P_{21}^{l-k} P_{21} P_{12} P_{12}^{l-k}\Omega_{12}^k f
    +    \sum_{k=0}^{l} \alpha_k^{(l)}
           \frac{ [n_2+l-2k] }{ [n_2+l-2k+1] }
           (12)^k P_{21}^{l-k} (12)
\Omega_{12} P_{12}^{l-k} \Omega_{12}^k f.    \]
By Proposition 3.7 it follows that
\[ f  =     \sum_{k=0}^{l} \alpha_k^{(l)}
(12)^k P_{21}^{l-k+1} P_{12}^{l-k+1}\Omega_{12}^k f
    +       \sum_{k=0}^{l} \alpha_k^{(l)}
(12)^{k+1} P_{21}^{l-k}
         \Omega_{12} P_{12}^{l-k} \Omega_{12}^k f
\frac{ [n_2+l-2k]   }{ [n_2+l-2k+1] } \times \]
\[  q^{l-k} \left(
    \frac{ [n_2+l-2k-1] }{ [n_2+l-2k  ] }
    \frac{ [n_2+l-2k-2] }{ [n_2+l-2k-1] }
        ...
    \frac{ [n_2-k]      }{ [n_2-k   +1] }
    \right)                        \]

\[ =       \sum_{k=0}^{l} \alpha_k^{(l)}
(12)^k P_{21}^{l-k+1} P_{12}^{l-k+1}\Omega_{12}^k f
 +       \sum_{k=0}^{l} \alpha_k^{(l)}
(12)^{k+1} P_{21}^{l-k}
         \Omega_{12} P_{12}^{l-k} \Omega_{12}^k f
\frac{ [n_2-k] }{ [n_2+l-2k+1] }               \]

\[ =       \sum_{k=0}^{l} \alpha_k^{(l)}
(12)^k P_{21}^{l-k+1} P_{12}^{l-k+1}\Omega_{12}^k f
   +       \sum_{k=0}^{l} \alpha_k^{(l)}
(12)^{k+1} P_{21}^{l-k}
       P_{12}^{l-k} \Omega_{12}^{k+1} f
\frac{ [n_2-k] }{ [n_2+l-2k+1] } \times            \]
\[  \left(
    \frac{ [n_2+l-2k-1] }{ [n_2+l-2k] }
    \frac{ [n_2+l-2k-2] }{ [n_2+l-2k-1] }
       ...
    \frac{ [n_2-   k]   }{ [n_2  - k+1] }
    \right)  \]

\[ =       \sum_{k=0}^{l} \alpha_k^{(l)}
(12)^k P_{21}^{l-k+1} P_{12}^{l-k+1}\Omega_{12}^k f \]
\[ +       \sum_{k=0}^{l} \alpha_k^{(l)}
(12)^{k+1} P_{21}^{l-k}
       P_{12}^{l-k} \Omega_{12}^{k+1} f
\frac{ [n_2-k]^2    }{ [n_2+l-2k] [n_2+l-2k+1] }          \]

\[ = \sum_{k=0}^{l}   ( \alpha_k^{(l)}
+ \alpha_{k-1}^{(l)} \frac{ [n_2-k+1]^2 }{
                            [n_2+l-2k+2] [n_2+l-2k+3]})
(12)^k P_{21}^{l+1-k} P_{12}^{l+1-k}\Omega_{12}^k f. \]
Therefore the statement is true for $l+1$ with (14).

The coefficients are uniquely determined by the relations (14).
Therefore it is sufficient to check the relations (14) for
\[\alpha_{k}^{(l)}: =
  \frac{\left[  \begin{array}{c} l \\ k \end{array} \right]
  \left[\begin{array}{c} n_2 \\ k \end{array} \right] }{
  \left[\begin{array}{c} l+n_2-k+1 \\ k \end{array} \right]}.\]
It follows from the identity
\[    [a][b-c]+[b][c-a]+[c][a-b] = 0  \]
with $a=l+1$, $b=k$, $c=l+n_2-k+2$ that
\[ [l+1][l+n_2-2k+2]= [k][n_2-k+1] + [l+n_2-k+2][l+1-k]. \]
Furthermore
\[   [l+1][l+n_2-2k+2]=  \]
\[   [l+n_2-k+2][l+1-k]
+ \frac{[k][n_2+l-2k+2][n_2+l-2k+3]}{[n_2-k+1] }
  \frac{[n_2-k+1]^2  }{[n_2+l-2k+2] [n_2+l-2k+3]}. \]
Therefore
{\tiny
\[ \frac{\left[  \begin{array}{c} l+1 \\ k \end{array} \right]
   \left[\begin{array}{c} n_2 \\ k \end{array} \right] }{
   \left[\begin{array}{c} l+n_2-k+2 \\ k \end{array} \right]} =
   \frac{\left[  \begin{array}{c} l \\ k \end{array} \right]
   \left[\begin{array}{c} n_2 \\ k \end{array} \right] }{
   \left[\begin{array}{c} l+n_2-k+1 \\ k \end{array} \right]} +
   \frac{\left[  \begin{array}{c} l \\ k-1 \end{array} \right]
   \left[\begin{array}{c} n_2 \\ k-1 \end{array} \right] }{
   \left[\begin{array}{c} l+n_2-k+1 \\ k-1 \end{array} \right]}
   (\frac{ [n_2-k+1]^2 }{ [n_2+l-2k+2] [n_2+l-2k+3] }). \]}
The assertion follows.$\bullet$

\begin{Rem} \rm
In the classical case
we can consider the vector space
generated by the monomials
$x_1^{m-i} y_1^i x_2^{n-j} y_2^j$,
$0\leq i\leq m$, $0\leq j \leq n$
as the representation
module $ V_m \otimes V_n $ of the Lie algebra $sl(2)$.
The elements
$ P_{12}^{n_1-k}
    \Omega_{12}^k   x_1^{m-i} y_1^i x_2^{n-j} y_2^j$,
$0 \leq k \leq {\rm min} (m,n)$
depend only on $x_2,y_2$.
One can show that the polars
\[ (12)^k  P_{21}^{n_1-k} P_{12}^{n_1-k}
    \Omega_{12}^k   x_1^{m-i} y_1^i x_2^{n-j} y_2^j  \]
span an irreducible subrepresentation $V_{m+n-2k}$.
Therefore we can use the quantised Gordan-Capelli series
for an explicite decomposition of the braided tensor product
of $V_m$ and $V_n$ (cf. \cite{Ma})
into irreducible representations of $U_q(sl(2,\RR)$, i.e.,
for $n=2$ the Gordan-Capelli series determines the structure
of the representation ring of $U_q (sl(2,\RR) $.
\end{Rem}

\begin{center}
\Large 4. $n$-FORMS
\end{center}
\setcounter{section}{4}
\setcounter{De}{0}
\setcounter{Ex}{0}
\setcounter{Rem}{0}
\setcounter{Sa}{0}

\begin{flushleft}
4.1 \it  Definition and Examples
\end{flushleft}

In this section we use the algebra $H_I$ in order
to introduce noncommutative binary homogeneous $n$-forms.

For the following we suppose for the algebra $H_I$ that
$0\in I$.
We will use the short notations
\[  x:=x_0,\ \ \ \ \  y:=y_0.   \]
Furthermore let $H_{I \backslash \{ 0 \} }$
be the subalgebra of $H_I$
which is generated by the elements
$x_i,y_i$,   $i\in I \backslash \{ 0 \}$.

\begin{De} {\rm
We say that the nonvanishing element
\begin{eqnarray} f:=                                        x^n        A_0 +
       \left[ \begin{array}{c} n \\ 1 \end{array}   \right] x^{n-1}y   A_1 +
       \left[ \begin{array}{c} n \\ 2 \end{array}   \right] x^{n-2}y^2 A_2 +
 ....  + y^n A_n
\end{eqnarray}
of $H_I$ with $A_i\in H_{I\backslash \{0\} }$
is a {\it n-form },
if $f\in H_I^{Inv} $. }
%We say that two $n$-forms $f$, $f'$
%are equivalent ($f \sim f' $)
%if there exists a $n$-form $f''$
%and elements $a,b\in H_{I \backslash \{ 0 \} }^{Inv}$
%with
%$ A_i = {A''}_i a$ and $ {A'}_i = {A''}_i b$ for $i=0,1,\cdots,n$.
\end{De}

\begin{Rem} \rm
The  representation of $f$  with coefficients
$A_i\in H_{I\backslash \{0\} }$ is unique.
This is a consequence of the PBW theorem.
We do not suppose anything about the commutation behaviour
of the elements $A_i$ (cf. Remark 4.2 below).
\end{Rem}

\begin{Ex} \rm
For the linear form
\[ f = x A + y B: = (0i) \]
($i\neq 0$) we find  $A=  q^{-\frac{1}{2}} y_i$,
          $B=- q^{ \frac{1}{2}} x_i$.
\end{Ex}

\begin{Ex} \rm
For the linear form
\[ f = x A + y B: =  (01)(23) -c (03)(12) \]
a calculation yields
\[ A =      c q^2  x_1 y_2 y_3 +
    (q^{-1}-c q^3) y_1 x_2 y_3 -
                   y_1 y_2 x_3, \]
\[ B=    -         x_1 x_2 y_3 +
    (q-c q^3)      x_1 y_2 x_3 +
     c q^{4}       y_1 x_2 x_3. \]
\end{Ex}

\begin{Ex} \rm
For the quadratic form
"with two different real zeros"
\[  f = x^2 A + x y [2]_{q^2} B + y^2 C := q(01)(02) \]
we find   $A= q y_1 y_2$,
          $B=\frac{-q^2 x_1 y_2-q  y_1 x_2}{q+q^{-1}}$
and       $C= q^3  x_1 x_2$,
(cf. Lemma 3.3 (i)).
(Later we will be able to present examples of quantised
real forms with complex zeros.)
\end{Ex}

\begin{Ex} \rm
For the quadratic form
"with two equal real zeros"
\[  f = x^2 A+[2]_{q^2} x y   B+y^2 C:= (01)^2 \]
we find $A=              y_1^2  $,
          $B= -      x_1   y_1    $
    and   $C=   q^2  x_1^2        $.
\end{Ex}

\begin{Ex} \rm
Consider $n$-forms of the type
\[ (0i_1) (0i_2) \cdots (0i_n) \]
which correspond to polynomials with real zeros.
Because of
$(0j)(0i)= q^{2} (0i)(0j)$, ($0<i<j$),
$(0j)(0i)= q^{4} (0i)(0j)$, ($i<0<j$) and
$(0j)(0i)= q^{2} (0i)(0j)$, ($i<j<0$)
a permutation of the indices
yields only a multiplication
of the form by a nonzero constant factor
(cf. Proposition 3.2 (i)).
Therefore it is sufficient to consider forms
$(0i_1) (0i_2) \cdots (0i_n)$ with $i_1\leq i_2\leq ... \leq i_n$.
We will endow each form of this kind
with a certain factor $q^c$, $c\in \NN $
such that the form becomes
a real (i.e. $*$-invariant) algebra element.
For example let $0 < i_1 < i_2 < \cdots < i_n$. Then
\[ f_{i_1,i_2,...,i_n}:=q^\frac{n(n-1)}{2} (0i_1)(0i_2)\cdots (0i_n) \]
is a real form.
\end{Ex}

\begin{Rem} \rm
The commutation behaviour of the coefficients of a $n$-form
depends on the concrete type of the form.
In Example 4.1 the coefficients $A,B$
of the linear form $f=xA+yB=(0i)$
have the commutation relations
\[           B A = q^2 AB.          \]
However for the linear form
$ f = x A + y B $ of Example 4.2 we have
\begin{eqnarray}
 BA=q AB+c q^\frac{7}{2} (q^2-1) (12)(13)(23).
\end{eqnarray}
One can show,
that there are no closed quadratic commutation relations
between $A$ and $B$.
\end{Rem}

%there are no linear combinations between
%$A^2, AB, AC, BA, B^2, BC, CA, CB, CA$ and therefore
%no quadratic relations between $A, B, C$.
%Higher order relations are
%
%\[  d^3  =  a   I_{ex}^2,               \ \ \ \ \
%    D_1  =  b   D_2       =  c D_3.  \]

\begin{Sa}
We have
\begin{eqnarray}
\begin{array}{ll}
L[A_i]=q^\frac{n}{2}    A_i, \ \ \ \ \ \ \   &
K[A_i]=q^\frac{n-2i}{2} A_i,                 \\
E[A_i]=-q^\frac{2i-n+2}{2} [i  ] A_{i-1},
\ \ \ \ \ \ \   &
F[A_i]=-q^\frac{-n-2  }{2} [n-i] A_{i+1}
\end{array}
\end{eqnarray}
($A_{-1}=A_{n+1}:=0$).
\end{Sa}

{\it Proof.} We have
\[   K(f) = \sum_{i=0}^n
\left[ \begin{array}{l} n \\ i \end{array} \right]
       q^\frac{-n+2i}{2} x^{n-i} y^i K( A_i ).  \]
By $K(f)=f$ and by the PBW theorem
the formula for $K$ follows. Furthermore
\[ E(f)=\sum_{i=0}^n
\left[ \begin{array}{l} n \\ i \end{array} \right]
  (         E(x^{n-i})      K(y^i    A_i)+
            K^{-1}(x^{n-i}) E(y^i) K(A_i)+
            K^{-1}(x^{n-i}    y^i) E(A_i)  ).       \]
By $E(x^m)=q^\frac{2-m}{2} [m] x^{m-1} y$ and
           $E(y^m)=0$ it follows that
\[= \sum_{i=0}^n
   \left[ \begin{array}{c} n \\ i \end{array} \right]
( q^\frac{2-n+i}{2} [n-i] x^{n-i-1} y q^\frac{i}{2}    y^i
                                      q^\frac{n-2i}{2} A_i+
  q^\frac{n-2i }{2} x^{n-i  }               y^i E(A_i)  ) \]
\[= \sum_{i=0}^n
\left[ \begin{array}{c} n \\ i-1 \end{array} \right]
     q                [n-i+1] x^{n-i}  y^{i}   A_{i-1} +
\left[ \begin{array}{c} n \\ i \end{array} \right]
     q^\frac{n-2i}{2}         x^{n-i}  y^i   E(A_i). \]
By $E(f)=0$ and by the PBW theorem
the formula for $E$ follows.
The proof for $F$ and $L$ is analogous.
$\bullet$

Therefore the action of $U_q(gl(2,\RR ))$ on the coefficients
is independent from the concrete realisation of the $n$-form.

\begin{flushleft}
4.2 \it  Geometric interpretation of $H_I$
\end{flushleft}

In Remark 2.3 we considered the generators $x_i,y_i$,
$i\in I$ as homogeneous coordinates
of points of the real line.
In this section we use these points in order to
construct further points.

\begin{De}{\rm We say that the ordered pair $(X,Y)$ of
           nonvanishing elements
           $X,Y \in H_{ I \backslash \{ 0 \} }^{Inv}$
           is a {\it point} if \\
(i)  $q^{-1/2} x_0 Y - q^{1/2} y_0 X\ \in H_I^{Inv}$ and \\
(ii) $X,Y$ are homogeneous with respect to all
     indices $i\in I$ of the same degree.  \\
We say, that the two points $(X_1,Y_1)$, $(X_2,Y_2)$
are equal if $X_1 Y_1^{-1} = X_2 Y_2^{-1}$ in $Q_I$.
}\end{De}

\begin{Rem} \rm
Le $q$ be not a root of unity.
Then $X,Y$ are elements of odd degree,
since $H_I^{Inv}$
is generated by the elements $(ij)$.
\end{Rem}

\begin{Ex} \rm
Pairs of generators $(x_i,y_i)$ are points. Furthermore let
$f=x_0 A+y_0 B$ be a linear form. Then
$(- q^{-1/2} B,  q^{1/2} A)$ is a point.
For the form
$f=  (01)(23) - c (03)(12)$
of Example 4.2 we have
\[ X=    q^\frac{-1}{2}                 x_1 x_2 y_3 -
        (q^\frac{ 1}{2}-c q^\frac{5}{2}) x_1 y_2 x_3 -
       c  q^\frac{7}{2}                 y_1 x_2 x_3,  \]

\[ Y=  c  q^\frac{5}{2}                 x_1 y_2 y_3 +
        (q^\frac{-1}{2}-c q^\frac{7}{2})  y_1 x_2 y_3 -
         q^\frac{1}{2}                 y_1 y_2 x_3. \]
\end{Ex}

\begin{Ex} \rm
Let $(X,Y)$ be a point and $a\in H_I^{Inv}$. Then we have
\[     (X,Y) = (Xa,Ya).  \]
\end{Ex}

%
%If we suppose further that $q$ is not a root of unity,
%we can represent $a$ by symbols.
%Since $X$ and $Y$ are homogeneous
%with respect to all indices of the same degree,
%it follows also
%\[     (X,Y) = (aX,aY)  \]
%by Lemma 3.1.
%

\begin{flushleft}
4.3 \it Polynomials, Elementary symmetric functions,
Newton relations
\end{flushleft}

In Remark 2.3 we have introduced the elements
$ v_i  = q^\frac{1}{2} x_i y_i^{-1}\in Q_I$
which obey the commutation relations
\[  v_j v_i = q^2 v_i v_j - ( 1 - q^2 )  v_i^2  \]
for $i<j$.

We define bracket symbols $[ij]$ by
\[       [ij]: =  v_i - v_j.       \]
They are related to the bracket symbols $(ij)$ by
\[       (ij)   =  y_i [ij] y_j.    \]
If we consider $v_0$ as an independent variable,
we can define
{\it polynomials with zeros}
$v_{i_1}$, $v_{i_2}$, ... , $v_{i_n}$ by
\[ p_{i_1,i_2,...,i_n}:=
          [0i_1][0i_2]\ \cdots \ [0i_n].  \]

\begin{Rem} \rm
The definition of $p_{i_1,i_2,...,i_n}$ depends on
the relative position of the brackets $[0i_j]$.
By Example 4.5
a permutation of the indices of the form
$(0i_1)(0i_2) \cdots (0i_n)$
yields only a multiplication by a constant factor.
We have not the same situation for  polynomials $p_{i_1,i_2,...,i_n}$.
For example we have
\[ p_{1,2} = [01][02]
           = q^{-2} v_0^2 -
             q^{-4} v_0 (q^2 v_1 + v_2) +
             q^{-4} v_1 v_2 \]
and
\[ p_{2,1} = [02][01]
           = q^{-2} v_0^2 -
             q^{-4} v_0 ( v_1 + q^2 v_2) +
             q^{-2} v_1 v_2 +
             q^{-4} (1-q^2)v_1^2. \]
\end{Rem}

\begin{Rem} \rm
The polynomial $p_{i_1,i_2,...,i_n}$
is related to the  form $f_{i_1,i_2,...,i_n}$ (cf. Example 4.5).
By Lemma 3.1 we have the commutation relations
\[  (0i) y_0^{-1} = \frac{1}{q}  y_0^{-1} (0i),\ \ \ \ \ \ \
    (0i) y_i^{-1} =          q   y_i^{-1} (0i),\ \ \ \ \  \ \
    (0i) y_j^{-1} =          q^3 y_j^{-1} (0i)          \]
for $0<i<j$.
It follows that
\[ p_{1,2,...,n} =          [01][02]\   \cdots \ [0n]
                 = y_0^{-1} (01)        y_1^{-1}
       y_0^{-1} (02) y_2^{-1} y_0^{-1} \ \cdots \
       y_0^{-1} (0n) y_n^{-1}                       \]
\[   = q^{n^2}
      y_0^{-1} y_1^{-1} y_0^{-1} y_2^{-1} \ \cdots \
      y_0^{-1} y_n^{-1} (01)(02)\ \cdots \ (0n)     \]
\[   = q^{\frac{n(3-n)}{2}}
      y_0^{-n} y_1^{-1} y_2^{-1}  \ \cdots \ y_n^{-1}
          f_{1,2,...,n}. \]
\end{Rem}

In the following we restrict the consideration to polynomials
$p_{1,2,...,n}$, $n\in\NN$.

\begin{Sa}
We have
\[ p_{1,2,...,n}= q^{n(1-n)} \sum_{j=0}^n   v_0^{n-j} A_j \]
with
\[ A_j = (-1)^j  \sum_{ 1\leq i_1 < i_2 < \ ... \ <i_j \leq n  }
                 q^{ 2j - 2(i_1+i_2+\ ... \ +i_j) }
      v_{i_1} v_{i_2}\ ... \ v_{i_j}          \]
for $j=1,..., n$ and $A_0=1$.
\end{Sa}
The statement can be shown by induction.

We say that $A_j$ is the {\it $j$-th symmetric function}
of $v_{i_1}$, $v_{i_2}$,\ ..., $v_{i_j}$.
Furthermore we define the {\it $j$-th  power sum} $P_j$, $j\geq 1$ by
\[  P_j =  v_1^j + q^{-2} v_2^j+ ... + q^{2-2n} v_n^j.    \]

\begin{Sa} Let $k\geq 1$ and $A_k:=0$ for $k>n$ or $k<0$.
Then we have
\[ P_k\ =\ \           -  q^{ 2 (k-1)(n-1) } (k)_{q^2} A_k\ -\
     \sum_{i=1}^{k-1} q^{ 2 (k-i)(n-1) }\     A_{k-i} P_{i}  \]
(Newton relations).
\end{Sa}

{\it Proof.} Let $n=1$. We have $P_k=v_1^k$, $A_0=1$, $A_1=-v_1$ and
$A_k=0$ for $k\neq 0,1$. It follows that
$ P_1 = - A_1 $
and
$ P_k = - A_1 P_{k-1} $
for $k\geq 2$.
Therefore the assertion is true for $n=1$.

We suppose that the assertion is true for $n$.
Let $A_i^{(n)}:=A_i$ be the coefficients of $p_{1,2,...,n}$.
By Proposition 4.2 and the definition of $P_j$ we have
\[   A_j^{(n+1)}  =  A_j^{(n)}  -  q^{-2n} A_{j-1}^{(n)}v_{n+1} \]
and
\[   P_j^{(n+1)}  =  P_j^{(n)}  +  q^{-2n} v_{n+1}^j
                  \ \ \ \ \  {\rm for}\ \   j\geq 1.\]
It follows that
\[ P_k^{(n+1)}      + q^{ 2(k-1)n } (k) A_k^{(n+1)} +
     \sum_{i=1}^{k-1} q^{ 2(k-i)n }     A_{k-i}^{(n+1)} P_{i}^{(n+1)}  \]

\[ = P_k^{(n)} + q^{-2n} v_{n+1}^k
 + q^{2(k-1)n}(k) A_k^{(n)} -  q^{2(k-1)n-2n} (k) A_{k-1}^{(n)} v_{n+1} \]
\[ +  \sum_{i=1}^{k-1}
     q^{ 2(k-i)n }     A_{k-i}^{(n)}            P_i^{(n)}
   -  \sum_{i=1}^{k-1}
     q^{ 2(k-i)n-2n }  A_{k-i-1}^{(n)}v_{n+1}   P_i^{(n)}    \]
\[ +  \sum_{i=1}^{k-1}
     q^{ 2(k-i)n-2n }  A_{k-i}^{(n)}            v_{n+1}^i
   -  \sum_{i=1}^{k-1}
     q^{ 2(k-i)n-4n }  A_{k-i-1}^{(n)}          v_{n+1}^{i+1}.   \]
The last row  gives
$q^{2(k-2)n} A_{k-1}^{(n)} v_{n+1} - q^{-2n} v_{n+1}^k$.
Because of the identity \\
$ v_{n+1} v_j^i = q^{2i} v_j^i v_{n+1} + (1-q^{2i}) v_j^{i+1} $
for $j<n+1$ we have
\[ = P_k^{(n)}
   + q^{2(k-1)n}(k) A_k^{(n)}
   + q^{2(k-2)n} ( 1-(k) ) A_{k-1}^{(n)} v_{n+1}
   + \sum_{i=1}^{k-1}
     q^{ 2(k-i)n }       A_{k-i  }^{(n)}    P_i^{(n)}  \]
\[ - \sum_{i=1}^{k-1}
     q^{ 2(k-i)n-2n+2i } A_{k-i-1}^{(n)}    P_i^{(n)} v_{n+1}
   - \sum_{i=1}^{k-1}
     q^{ 2(k-i)n-2n } (1-q^{2i}) A_{k-i-1}^{(n)}    P_{i+1}^{(n)}. \]
%Since
%$  \sum_{i=1}^{k-1}
%  q^{ 2(k-i)n-2n } (1-q^{2i  }) A_{k-i-1}^{(n)} P_{i+1}^{(n)}=
%   \sum_{i=2}^{k}
%  q^{ 2(k-i)n    } (1-q^{2i-2}) A_{k-i  }^{(n)} P_{i  }^{(n)}$
%we have
We unite the first and the third sum and rewrite the second sum.
\[ = P_k^{(n)}
               + q^{2(k-1)n}(k) A_k^{(n)}
               + q^{2(k-2)n} ( 1-(k) ) A_{k-1}^{(n)} v_{n+1}
 + \sum_{i=1}^{k-1}
     q^{ 2(k-i)n+2i-2 }  A_{k-i}^{(n)}  P_i^{(n)}          \]
\[   - (1-q^{2k-2})        A_{0  }^{(n)}  P_k^{(n)}
     -q^{2k-2} \sum_{i=1}^{k-2}
     q^{ 2(k-1-i)(n-1)} A_{k-i-1}^{(n)}    P_i^{(n)} v_{n+1}
       -q^{2k-2} P_{k-1}^{(n)}.                            \]
Because the assertion is true for $n$
we can replace the first sum by
$ -q^{2k-2} (P_k^{(n)} + q^{2(k-1)(n-1)} A_k^{(n)})$
and the last two terms by
$(k-1)_{q^2} q^{2(k-2)(n-1)+2k-2} A_{k-1}^{(n)} v_{n+1}$.
We obtain
\[ = P_k^{(n)}
               + q^{2(k-1)n}(k)_{q^2} A_k^{(n)}
               + q^{2(k-2)n} ( 1-(k)_{q^2} ) A_{k-1}^{(n)} v_{n+1}
               - q^{2k-2} (P_k^{(n)} + q^{2(k-1)(n-1)} A_k^{(n)})       \]
\[   - (1-q^{2k-2})        A_{0  }^{(n)}  P_k^{(n)}
        +  (k-1)_{q^2} q^{2(k-2)(n-1)+2k-2} A_{k-1}^{(n)} v_{n+1}   \]
\[ =  (q^{2(k-2)n} ( 1-(k)_{q^2} )
     + (k-1)_{q^2} q^{2(k-2)(n-1)+2k-2}) A_{k-1} v_{n+1} = 0. \]
The assertion follows.$\bullet$

\begin{center}
\Large 5. INVARIANTS OF FORMS
\end{center}
\setcounter{section}{5}
\setcounter{De}{0}
\setcounter{Ex}{0}
\setcounter{Rem}{0}
\setcounter{Sa}{0}

\begin{flushleft}
5.1 \it Definition
\end{flushleft}

We consider a $n$-form $f=\sum_{i=0}^n
\left[ \begin{array}{l} n \\ i \end{array} \right]
x^{n-i} y^i A_i $.

\begin{De}  {\rm
We say that the polynomial expression of the coefficients of $f$
\[ I_f = I_f (A_0,...,A_n) \]
is an {\it invariant},
if $I_f \in H_{I \backslash \{ 0 \} }^{Inv}$.
We say that an invariant $I_f$
is an {\it invariant of degree $k$},
if $I_f$ is a polynomial of degree $k$
in the coefficients $A_i$.
We say that the polynomial expression
\[ C_f = C_f ( A_0,...,A_n ,x_0,y_0 ) \]
is a {\it covariant} of $f$ of degree $m$,
if $C_f$ is a $m$-form.
%if  $L(A_i)= q^\frac{m }{2} A_i$ with $m\in \ZZ$
%for $i=0,...,n$
%and $L(I_f)= q^\frac{km}{2} I_f$.
}
\end{De}

\begin{Ex} \rm
Let $f= q (01)(02)$ (cf. Example 4.3) and
\[ I_f :=
\left(
\frac{1}{q} A C+\frac{1}{q^3} C A \right)
- \frac{1}{q} [2] B^2. \]
A calculation yields
\[ I_f = - \frac{1}{ [2] } (12)^2.  \]
Therefore $I_f$ is an invariant of degree 2.
\end{Ex}

Sums and scalar multiples of invariants of degree $k$
are again invariants of degree k.
By (2) the product of two invariants
of the degrees $k$ and $l$ is an invariant of degree $k+l$.
Therefore the invariants of a $n$-form $f$ form a
graded algebra.

Furthermore we define invariants of a collection of forms.

\begin{De}
We say that a polynomial expression
$I_{f',f'',...}$
of the coefficients
${A'}_i, {A''}_i, ... $  of the binary forms
$f', f'', ...  \in H_I$ of the degrees $n'$, $n''$, ...
is a {\it simultaneous invariant},
if $I_{f',f'',...} \in H_{I \backslash \{ 0 \} }^{Inv}$.
Analogous we define common covariants.
\end{De}

\begin{Ex} \rm
Let $f'=(01)$, $f''=(02)$
(cf. Example 4.1) and
\[   I_{f',f''} =
     q^{-\frac{1}{2}} B' A'' - q^\frac{1}{2} A' B''. \]
We have
\[   I_{f',f''} =  -  (12)
     \in H_{I \backslash \{ 0 \} }^{Inv}.  \]
\end{Ex}

\begin{flushleft}
5.2 \it  Universal invariants
\end{flushleft}

The Remark 4.2 demonstrates
us that the structure of the coefficient algebra depends
on the concrete type of the $n$-form. Therefore it is useful to define the
notions of {\it universal $n$-forms} and {\it universal invariants}.

\begin{De} \rm
We say that the expression
\[ {\cal F}:=                                      x^n        {\cal A}_0 +
\left[ \begin{array}{c} n \\ 1 \end{array} \right] x^{n-1}y   {\cal A}_1 +
\left[ \begin{array}{c} n \\ 2 \end{array} \right] x^{n-2}y^2 {\cal A}_2 +
                                      ....       + y^n        {\cal A}_n \]
is an {\it universal $n$-form}, if the ${\cal A}_i$ are abstract variables.
We obtain a $n$-form as a {\it realisation} from an universal form,
if we insert the elements $A_i \in H_I$ for ${\cal A}_i$.
\end{De}

Consider $k$ universal $n$-forms  ${\cal F}',...,{\cal F}^{(k)}$.
We form the free algebra $F_I$
generated by their coefficients ${\cal A}^{(j)}_i$.
For the realisation $f',...,f^{(k)}$
we have the homomorphism $h_{f',...,f^{(k)} }: F_I \rightarrow H_I$
determined
by $h_{f',...,f^{(k)} } ( {\cal A}^{(j)}_i ):= A^{(j)}_i$.

The formulas (17) determine also an $U_q(gl(2,\RR))$-action on $F_I$
if we replace $A_i$ by ${\cal A}_i$. We have
\begin{eqnarray}
X(h_{f',...,f^{(k)} } (a))=h_{f',...,f^{(k)} } (X(a)), \ \ \ \ \
\end{eqnarray}
for $X \in  U_q(gl(2,\RR )$ and $a\in F_I$.

\begin{De} \rm
We say that ${\cal I}\in F_I$ is an
{\it universal invariant} of ${\cal F}$ if ${\cal I}$ is an
invariant element of the module algebra $F_I$
(i.e. $E({\cal I})=F({\cal I})=0$ and $K({\cal I})={\cal I}$).
Analogous we define {\it universal simultaneous invariants}
and {\it universal covariants}.
\end{De}

\begin{Ex} \rm
Consider the universal form
\[ x^2 {\cal A} + [2]_{q^2} x y {\cal B} + y^2 {\cal C} \]
and
\[{\cal I} := \frac{1}{q} {\cal A C}+
        \frac{1}{q^3}    {\cal C A}
      - \frac{1}{q}  [2] {\cal B}^2. \]
By (17) and (18) we have
$E({\cal A})=0,
 E({\cal B})=-q {\cal A},
 E({\cal C})=-q^2 {\cal B}$ and  \\
$K({\cal A})=q      {\cal A},
 K({\cal B})=       {\cal B},
 K({\cal C})=q^{-1} {\cal C}$.
Therefore we obtain $K({\cal I})={\cal I}$ and
\[ E({\cal I})  = \frac{1}{q}      K^{-1} ({\cal A}) E({\cal C})+
                  \frac{1}{q}      E      ({\cal A}) K({\cal C})+
                  \frac{1}{q^3}    K^{-1} ({\cal C}) E({\cal A})  \]
\[              + \frac{1}{q^3}    E      ({\cal C}) K({\cal A})-
                  \frac{1}{q}  [2] K^{-1} ({\cal B}) E({\cal B})-
                  \frac{1}{q}  [2] E      ({\cal B}) K({\cal B})  \]
\[ = -[2] {\cal AB}-[2] {\cal BA}+[2] {\cal AB}+[2] {\cal BA}=0 . \]
Similarly we find  $F({\cal I})=0$.
\end{Ex}

\begin{Le}
     Let ${\cal I}$ be an universal invariant of $\cal F$
     and let $f$ be a realisation of $\cal F$.
     Then $I_f := h_f ({\cal I})$ is an invariant of $f$.
\end{Le}

{\it Proof.}   The assertion follows from (18).
$\bullet $

\begin{Rem} \rm
Let $\cal I$ be a simultaneous universal invariant
of the universal n-forms
$\cal F',F'',...$ of the same degree $k$.
If we identify the coefficients of the different forms
with the coefficients of $\cal F$,
then $\cal I$ becomes an invariant of the form $\cal F$.
This follows, since the form of the
$U_q(sl_2)$-action on the coefficients of the
different forms ${\cal F}^{(j)}$ is independent on $j$.
\end{Rem}

\begin{flushleft}
5.3 \it  The symbolic method of Clebsch and Gordan
\end{flushleft}

Consider the universal $n$-form
\[ {\cal F} =                                      x^n        {\cal A}_0 +
\left[ \begin{array}{c} n \\ 1 \end{array} \right] x^{n-1}y   {\cal A}_1 +
\left[ \begin{array}{c} n \\ 2 \end{array} \right] x^{n-2}y^2 {\cal A}_2 +
                                      ....       + y^n        {\cal A}_n. \]
Our aim is to construct
universal invariants ${\cal I}$ of degree $k$.

We consider the collection of $k$ realisations of $\cal F$
\[ f':=(01)^n,\ \ f'':=(02)^n,\ \  ... ,\ \  f^{(k)}:=(0k)^n. \]
By  Lemma 3.3 (ii)
the forms have the coefficients
\begin{eqnarray}
{ A^{(j)} }_i = (-1)^i
q^{ \frac{1}{2} n^2 + (i-n)(1+i) } x_j^i y_j^{n-i}.
\end{eqnarray}

Consider an element $d$ of $H^{Inv}_I$
which is homogeneous of degree $n$
with respect to $1,2,...,k$.
(For example $(12)(13)(23)$ is
homogeneous of degree 2 with respect to 1, 2, 3.)

Because of the PBW theorem for $H_I$
and because the application of the reduction rules (5)
does not change the degree of
homogenity with respect to $1,2,...,k$
we can express $d$ as a linear combination of monomials
\begin{eqnarray}
   x_1^{ n- i_1 } y_1^{ i_1 }
   x_2^{ n- i_2 } y_2^{ i_2 } ...
   x_k^{ n- i_k } y_k^{ i_k }.
\end{eqnarray}
Using (19) we can represent $d$ in terms of the
coefficients ${A^{(j)}}_i $. We  obtain the expression
\begin{eqnarray*}
d = \sum c_{i_1,i_2,...,i_k}
A^{(1)}_{i_1} A^{(2)}_{i_2} ... A^{(k)}_{i_k}.
\end{eqnarray*}
Furthermore we consider the elements of $F_I$
\[    {\cal I}' = \sum c_{i_1,i_2,...,i_k}
{\cal A}^{(1)}_{i_1} {\cal A}^{(2)}_{i_2} ... {\cal A}^{(k)}_{i_k} \]
and
\[ E({\cal I}') =\sum  d_{i_1,i_2,...,i_k}
{\cal A}^{(1)}_{i_1} {\cal A}^{(2)}_{i_2} ... {\cal A}^{(k)}_{i_k}. \]
By (18) it follows that
\[   0 = E(d) = E(h_{f',...,f^{(k)} } ({\cal I}')) =
                h_{f',...,f^{(k)} } (E({\cal I}')) =
     \sum d_{i_1,i_2,...,i_k}
     A^{(1)}_{i_1} A^{(2)}_{i_2} ... A^{(k)}_{i_k}    \]
\[ =\sum d_{i_1,i_2,...,i_k}
         q^{ \frac{1}{2} k_{i_1,i_2,...,i_k} }
            x_1^{ n- i_1 } y_1^{ i_1 }
            x_2^{ n- i_2 } y_2^{ i_2 } ...
            x_k^{ n- i_k } y_k^{ i_k }     \]
with certain numbers $k_{i_1,i_2,...,i_k}\in \ZZ$. We have $|
q^{\frac{1}{2} k_{i_1,i_2,...,i_k} } | = 1$. By the PBW theorem we
have $d_{i_1,i_2,...,i_k}=0$, i.e. $E( {\cal I}' )=0$. Similarly
we encounter $F( {\cal I}' )=0$ and $K( {\cal I}' )= {\cal I}' $.
Therefore $ {\cal I}'$ is an universal simultaneous invariant of
$k$ universal $n$-forms ${\cal F}'$, ${\cal F}''$, ... , ${\cal
F}^k$.

By Remark 5.1
we obtain again an invariant if we identify
${\cal A}^{(j)}_i$ with ${\cal A}_i$.
Therefore
\[ {\cal I} = \sum c_{i_1,i_2,...,i_k}
   {\cal A}_{i_1} {\cal A}_{i_2} ... {\cal A}_{i_k}  \]
is an universal invariant of degree $k$
of the universal $n$-form $\cal F$.

Therefore we have demonstrated the following Theorem.

\begin{Th}
Let $\cal F$ be an universal $n$-form
and let $d$ be an element of $H_I^{Inv}$
which is homogeneous of degree $n$ with respect to $1,2,...,k$.
Then the above construction gives
an universal invariant of $\cal F$ of
degree $k$.
\end{Th}

Analogous we can use the Clebsch-Gordan method in order to construct
universal covariants, universal common invariants
and universal common covariants.

We say that the element of $H^{Inv}_I$
\[ d=   (i_1 i_2)(i_3 i_4) ... (i_{kn-1} i_{kn})   \]
is a {\it symbol} of $\cal F$, if $d$
is homogeneous of degree $n$
with respect to $1,2,...,k$.

The simplest symbols are
$(12)^2$, $(12)(13)(23)$
for the quadratic form,
$(12)^3$
for the cubic form and
$(12)^4$, $(12)^2 (13)^2 (23)^2$
for the quartic form.

If $q^{2d}\neq 1$ for all $d\in \NN$, then we can represent
the generating element $d$ by a symbol
by Theorem 3.6.

Using Lemma 5.1, we obtain invariants of degree $k$ for
the realisations of the universal $n$-form.

\begin{Ex} \rm
Consider the universal quadratic form
\[ f = x^2 {\cal A} +[2]  x y {\cal B} +  y^2 {\cal C}.  \]
For the invariant $(ij)^2$ we have (cf. Lemma 3.3 (ii))
\[       (ij)^2=                 x_i^2 y_j^2
                       -    [2]  x_i y_i x_j y_j
                       +    q^2  y_i^2 x_j^2.      \]
Consider the realisations  $f=(0i)^2$, $i=1,2$  with
\[ A^{(i)} =      y_i^2,   \ \ \ \ \
   B^{(i)} = -    x_i y_i, \ \ \ \ \
   C^{(i)} = q^2  x_i^2.              \]
We obtain for the symbol  $d:=1/q (12)^2$
\[    \frac{1}{q}   (12)^2
=\frac{1}{q} x_1^2 y_2^2-\frac{1}{q} [2] x_1 y_1 x_2 y_2+q y_1^2 x_2^2 \]
\[  = \frac{1}{q^3} C' A'' -\frac{[2]}{q } B' B''+ \frac{1}{q} A' C''. \]
Therefore we have an universal common invariant
of ${\cal F}'$ and ${\cal F}''$.
\[ {\cal I}' :=
            \frac{1}{q^3}  {\cal C}' {\cal A}''
       - (1+\frac{1}{q^2}) {\cal B}' {\cal B}''
       +    \frac{1}{q}    {\cal A}' {\cal C}''. \]
If we identify the coefficients of ${\cal F}'$ and ${\cal F}''$
we obtain again the universal invariant
of Example 5.3.
\[                       {\cal I}
      =   \frac{1}{q^3}  {\cal C A }
        - \frac{[2]}{q}  {\cal B B }
        + \frac{1}{q}    {\cal A C }.  \]
\vspace{0.5cm}
\end{Ex}

We can obtain all universal invariants
by the symbolic method.

\begin{Sa}
Every universal invariant of the $n$-form $\cal F$
is a linear combination of invariants,
which are derived from elements $d \in H_I^{Inv}$
by the above method. If $q^{2d}\neq 1$ for all $d\in \NN$,
then we can derive all universal invariants
from symbols of $\cal F$.
\end{Sa}

{\it Proof.}
Let
\[   {\cal I} = \sum c_{i_1,i_2,...,i_k}
{\cal A}_{i_1} {\cal A}_{i_2} ... {\cal A}_{i_k}  \]
be an universal invariant of degree $m$ of
the $n$-form $\cal F$
and let
\[ f^{(j)}=(0j)^n =
\sum_{i=1}^n
\left[ \begin{array}{c} n \\ i \end{array} \right]
x^{n-i} y^i A_i^{(j)},\]
with $j=1,...,k$ be realisations of $\cal F$.
We set
\[   d := \sum c_{i_1,i_2,...,i_k}
{A}_{i_1}^{(1)} { A}_{i_2}^{(2)} ... { A}_{i_k}^{(k)}.\]
$d=h_{ f',...,f^{(k)} } ({\cal I})$
is an invariant element of  $H_{ I \backslash \{ 0 \} }$.
From the construction of $d$ it follows that
we obtain again ${\cal I}$ from $d$
by the symbolic method.
Using Theorem 3.6
we can express $d$ as a linear combination
of symbols of $\cal F$ if $q^{2d}\neq 1$ for all $d\in \NN$.
$\bullet$

\begin{Rem} \rm
We can modify the symbolic method, if we use the elements
\[
   x_{\pi (1)}^{ n- i_1 } y_{\pi (1)}^{ i_1 }
   x_{\pi (2)}^{ n- i_2 } y_{\pi (2)}^{ i_2 } ...
   x_{\pi (k)}^{ n- i_k } y_{\pi (k)}^{ i_k },
\]
where $\pi\in S_k$ is a permutation, as PBW basis elements
instead of (20).
It is a consequence of Proposition 5.3
that we obtain no further invariants by this modification.
\end{Rem}

\begin{flushleft}
5.4 \it The case of $q$ being a root of unity
\end{flushleft}

In this section we give some remarks for the case
of $q$ being a root of unity. We mention that the
Clebsch-Gordan symbolic method gives all universal
invariants even in this case.
However, if $q^{2d}=1$ for $d\in \NN$,
we have invariant elements of $H_I$, which cannot
be represented by symbols.

\begin{Sa}
Let $q$ be a primitive $2d$-$th$ root of unity
($d\geq 2$).
Then $H_{ \{ 1 \} }^{Inv}$ is generated by the elements
$      x_1^{k d} y_1^{k' d} $
with $ k, k' \in \NN$ and $4|(k-k')$.
\end{Sa}

{\it Proof.}
Let   $ f = \sum_{i,j} a_{ij} x_1^{i} y_1^{j}$
be an invariant element. We have
\[Kf =\sum_{i,j}a_{ij} q^\frac{ -i+j}{2} x_1^{i} y_1^{j} = f,         \]
\[Ef =\sum_{i,j}a_{ij} [i] q^\frac{2-i+j}{2} x_1^{i-1} y_1^{j+1} = 0  \]
and
\[ Ff = \sum_{i,j}a_{ij} [j] q^\frac{ i-j}{2} x_1^{i+1} y_1^{j-1} = 0.\]
Note that $[i]=0$, if and only if $d|i$.
Therefore  $H_{ \{ 1 \} }^{Inv}$
is spanned by the monomials
$x^{k d} y^{k' d}$ with $4| (k-k')$.$\bullet$

Now we can use the Clebsch-Gordan symbolic method
in order to determine the universal invariants of degree 1.
Let again $q$ be a primitive $2d$-$th$ root of unity
($d\geq 2$),
let $\cal F$ be  the universal $n$-form and let
\[ f=(01)^n=\sum_{k=0}^n
     \left[  \begin{array}{c} n \\ k \end{array}   \right]
     (-1)^k     q^{\frac{1}{2} n^2 +(k-n)(1+k)}
              x_0^{n-k} y_0^k x_1^{k} y_1^{n-k} \]
be a realisation.
For the invariant element $x_1^{k d} y_1^{k' d}$
with $4|k-k'$ we have $(k+k')d=n$.

Therefore ${\cal F}$ has universal invariants only for
even $\frac{n}{d}$
and by the Clebsch-Gordan symbolic method certain coefficients
are obtained as universal invariants:

\noindent
Case 1.: $\frac{n}{d} \equiv  2\ mod\ 4$:
\begin{eqnarray}
 {\cal A}_d,   \ \          {\cal A}_{3d},\ \
 {\cal A}_{5d},\ \cdots \ , {\cal A}_{ n-d },
\end{eqnarray}
Case 2.: $\frac{n}{d} \equiv  0\ mod\ 4$:
\begin{eqnarray}
 {\cal A}_0,   \ \          {\cal A}_{2d},\ \
 {\cal A}_{4d},\ \cdots \ , {\cal A}_{ n }.
\end{eqnarray}
By Proposition 5.4 these
are all universal invariants of degree 1.
We note that there are no invariants of degree 1
for $q^{2d} \neq 1$  and for $q^2=1$.

\begin{flushleft}
5.5 \it  Polarisation
\end{flushleft}

Let $f$ be the $n$-form
\begin{eqnarray}
f=\sum_{i=0}^n
\left[ \begin{array}{l} n \\ i \end{array} \right]
x^{n-i} y^i A_i
\end{eqnarray}
with $A_i \in H_{ I\backslash \{ 0,p \} }$
and $k\in \NN$.
We define the $k$-$th$ {\it polar} $f^{(k)}$
of $f$ by
\[  f^{(k)}:=  (P_{0p})^k f.   \]
From this definition we infer the properties
$f^{(n)}\in H_{ I\backslash \{ 0 \} }$,
$f^{(k)}=0$ for $k>n$ and
\[(f^{(k)})^{(l)} = (f^{(l)})^{(k)} = f^{(k+l)} \]
for $k,l\in \NN$.

\begin{Sa}
Let $f$ be the $n$-form (23).
%with $A_i \in H_{ I\backslash \{ 0,p \} }$.
We can represent $f$ in the form
\[    f=\sum_{i=0}^n
\left[ \begin{array}{l} n \\ i \end{array} \right]
A_i^- x^{n-i} y^i A_i^+    \]
with $A_i^- \in H_{ \{ i<0 \} }$
and  $A_i^+ \in H_{ \{ i>0 \} }$.
%Further let
%\[  (02)^n = \sum_{i=0}^n x_0^{n-i} y_0^{i} A_i    \]
%a realisation of $\cal F$.
Furthermore let
$C_{i}^{n,k} \in H_{\{ 0,1\} }$ be
the elements which are uniquely determined
by the PBW expansion
\[
\sum_{i=0}^n  C_{i}^{n,k} A_i :=
(02)^{n-k} (12)^k .
\]
with $A_i\in H_{ \{ 2 \}  }$.
Then we find  for the $k$-th polar of $ f$
\[{ f}^{(k)} =
\sum_{i=0}^n A_i^- C_{i}^{n,k} A_i^+. \]
\end{Sa}

\begin{Ex} \rm
Let $p=1$. For the quadratic form
\[  { f} =  x_0^2   {  A} +  [2]  x_0 y_0 {  B} + y_0^2  {  C}  \]
with $A,B,C \in H_{ \{ i\in I|i>1 \} }$,
Lemma 3.3(i) implies that
the first polar is
\[   { f}^{(1)}= x_0 (x_1 { A} +  q y_1 { B})
                   + y_0 (x_1 { B} +    y_1 { C}).\]
For $f=q(02)(03)$
we obtain the linear form
\[   f^{(1)}=
        q       x_0 x_1    y_2 y_3
      - \frac{1}{[2]}
       (q^2     x_0 y_1 (q x_2 y_3 + y_2 x_3)
      + q       y_0 x_1 (q x_2 y_3 + y_2 x_3))
      + q^3     y_0 y_1    x_2 x_3             \]
\[   = \frac{1}{[2]}((02)(13) + (12)(03)).     \]
\end{Ex}

\begin{Ex} \rm
 Let $p=1$.
For the cubic form
\[ { f}
   =     x_0^3   { A} + [3] x_0^2 y_0 { B} +
     [3] x_0 y_0^2 { C} + y_0^3       { D} \]
with $A,B,C,D \in H_{ \{ i\in I|i>1 \} }$
we find
\[ { f}^{(1)}
   =      x_0^2 x_1                         (  { A}
   +    (q^2 x_0^2 y_1+(q^{-1}+q)x_0 y_0 x_1)  { B}
   +    ((1+q^2)x_0 y_0 y_1+y_0^2 x_1 )        { C}
   +      y_0^2 y_1                            { D} \]
\[ =     x_0^2       ( x_1 { A} +  q^2 y_1 { B} )
       + x_0 y_0 [2] ( x_1 { B} +  q   y_1 { C} )
       + y_0^2       ( x_1 { C} +      y_1 { D} )
\]
and
\[ { f}^{(2)}
   =  x_0 x_1^2                               { A}
   + ((1+q^2) x_0 x_1 y_1      + y_0 x_1^2)   { B}
   + (q^2 x_0 y_1^2 +(q^{-1}+q)y_0 y_1 x_1)   { C}
   +  y_0 y_1^2                               { D}
\]
\[
 = x_0 (x_1^2 { A}+q[2] x_1 y_1 { B}+q^2 y_1^2 { C})
 + y_0 (x_1^2 { B}+ [2] y_1 x_1 { C} +   y_1^2 { D}).\]
For $f=q^3 (01)(02)(03)$ we obtain the quadratic form
\[ f^{(1)}=\frac{q}{[3]}((02)(03)(14)+(02)(13)(04)+(12)(03)(04)) \]
and the linear form
\[ f^{(2)}=\frac{q}{[3]}((02)(13)(14)+(12)(03)(14)+(12)(13)(04)).\]
\end{Ex}

\begin{center}
\Large 6. APPLICATION OF THE SYMBOLIC METHOD
\end{center}
\setcounter{section}{6}
\setcounter{De}{0}
\setcounter{Ex}{0}
\setcounter{Rem}{0}
\setcounter{Sa}{0}

\begin{flushleft}
6.1 \it  Linear Forms
\end{flushleft}

We consider the universal linear form
\[  {\cal F} =  x {\cal A} + y {\cal B}.     \]
In the classical case
the linear form $\cal F$ has no invariants.
For the realisations
$f = x A + y B: = (0i)$
($i=1,2$) we have $A=  q^{-\frac{1}{2}} y_i$ and
                  $B=- q^{ \frac{1}{2}} x_i$.
Therefore, for the symbol
$(12)=q^{-\frac{1}{2}} x_1 y_2 - q^{\frac{1}{2}} y_1 x_2 $
we obtain the universal invariant
\[  {\cal I}_1 =
- q^{-\frac{1}{2}}{\cal B A} + q^\frac{1}{2} {\cal A B}. \]

\begin{Ex} \rm
 Let $q^{2d}\neq 1$ for all $d\in \NN$.
We consider the simplest realisation
$f=(0i)$, $i\neq 0$ (cf. Example 4.1).
For an arbitrary universal invariant $\cal I$
we obtain the realisation $I_f =0$.

{\it Proof.}
Since each universal invariant $\cal I$
is generated by a symbol,
$\cal I$ is homogeneous
with respect to $\cal A$ and $\cal B$
of the same degree.
Therefore we find
$I_f= c\ x_i^k y_i^k$, $k\geq 1$ with $c\in \CC$.
The assertion follows from the formulas (3).
$\bullet$
\end{Ex}

The next Example demonstrates, that there are nontrivial
realisations of universal invariants for linear forms.

\begin{Ex} \rm
For the realisation
                $f =  (01)(23) - c (03)(12) $
                (cf. Example 4.2)
                we obtain
\[     I_{1,f} = c(q^3-q^5) (12)(13)(23) \neq 0,   \]
(cf. formula (16)). One can show, that all invariants of
$f$ have the form $\alpha I_{1,f}^k$
with $\alpha \in \CC$ and $k\in \NN$.
The proof is similar to the proof of Proposition 6.3
(cf. below).
\end{Ex}

A main task of the classical and the modern
representation theory is to find finite generating sets
of the algebras of invariants of forms.
The following proposition shows,
that there is no finite generating set
for universal invariants
even in the case of the universal linear form.
In this article we give finitness theorems
for realisations of universal forms only for
simple Examples (cf. the Examples 6.1 and 6.2 above
and the Propositions 6.2, 6.3 below).

\begin{Sa}
Let $q^{2d}\neq 1$ for all $d\in \NN$.
The algebra of universal invariants of the universal linear form
has not a finite set of generating elements.
\end{Sa}

{\it Proof.}
Consider the sequence of symbols
\[  (1,2n)(2,2n-1)(3,2n-2) \cdots (n,n+1)  \]
with $n=1,\ 2,\ \ \cdots \ \ \ $.
By Lemma 3.3 (iv) we obtain by the symbolic
method the sequence of universal invariants
\[ {\cal I}_n =
\sum_{ i_1,i_2, \cdots , i_n = 0 }^1
 - q^{ \frac{3}{2} n(n-1) }  (-q)^{ i_1+i_2+ \cdots +i_n }
{\cal U}^{(i_1)}
{\cal U}^{(i_2)}
\cdots
{\cal U}^{(i_n)}
{\cal U}^{(1-i_n)}
{\cal U}^{(1-i_{n-1})}
\cdots
{\cal U}^{(1-i_1)}
\]
of degree $2n$
with ${\cal U}^{(0)}:= {\cal B}$
and  ${\cal U}^{(1)}:= {\cal A}$.
The ${\cal I}_n$ are elements of the free algebra $F_I$.
Since each  universal invariant $\cal I$
is generated by a symbol,
$\cal I$ is a linear combination of monomials,
which are homogeneous
with respect to $\cal A$ and $\cal B$
of degree $k\in \NN$.
Because ${\cal I}_n$ contains the monomial
$- q^{ \frac{3}{2} n(n-1) } {\cal B}^n {\cal A}^n$,
it is impossible to express ${\cal I}_n$ in terms
of universal invariants ${\cal I}$
of degree $2k$ with $k<n$.
The assertion follows $\bullet$.

\begin{flushleft}
6.2 \it Quadratic forms
\end{flushleft}

In Example 5.4 we obtained for the symbol  $d := 1/q (12)^2$
by the symbolic method the universal invariant
\[       {\cal I}_{1}
      =   \frac{1}{q^3}  {\cal C A }
        - \frac{[2]}{q}  {\cal B B }
        + \frac{1}{q}    {\cal A C }.  \]

\begin{Ex} \rm
In Example 5.1 we found for the realisation $ f=   q (01)(02) $
\begin{eqnarray}
   I_{1,f} =   - \frac{1}{[2]} (12)^2.
\end{eqnarray}
For $ f= (01)^2 $ we find  $I_{1,f}=0$.
\end{Ex}

Now we consider an invariant without a classical analog.
Let ${\cal I}_{2}$ be the invariant, which  arises from the symbol
$-\frac{1}{q^4} (12)(13)(23)$. We obtain
\[{\cal I}_2 =
     q^\frac{  -5}{2}{\cal A} {\cal B} {\cal C}
 -   q^\frac{  -9}{2}{\cal A} {\cal C} {\cal B}
 -   q^\frac{  -9}{2}{\cal B} {\cal A} {\cal C}
 +   q^\frac{  -9}{2}{\cal B} {\cal C} {\cal A}
 +   q^\frac{  -9}{2}{\cal C} {\cal A} {\cal B}
 -   q^\frac{ -13}{2}{\cal C} {\cal B} {\cal A} \]
\[-(q^\frac{-7}{2}-q^\frac{ -15}{2}){\cal B}{\cal B} {\cal B}.\]

\begin{Rem} \rm
${\cal I}_2$ is
associated to the combinant of three quadratic forms
(cf. below).
\end{Rem}

For the realisation $f=q(01)(02)$ we have
\begin{eqnarray}
I_{2,f} = ( q^2 - \frac{1}{q^2} ) (12)^3.
\end{eqnarray}
We note that $I_{2,f} \neq 0$ if $q^4 \neq 1$.

Now we consider invariants of degree four.
The fourth order symbols are ge\-ne\-ra\-ted
by $Z_1^2, Z_1 Z_3, Z_3^2$
(c.f. Remark 3.1). We obtain by the symbolic method
\[   {\cal I}_{3,1}=   q^8    {\cal ACAC}
            -q^{10} [2]       {\cal ACBB}
            + q^6             {\cal ACCA}
            -q^{10} [2]       {\cal BBAC}
           + q^{12} [2]^2     {\cal BBBB}   \]
\[         - q^8 [2]          {\cal BBCA}
           + q^6              {\cal CAAC }
           - q^8 [2]          {\cal CABB }
           + q^4              {\cal CACA }, \]

\[ {\cal I}_{3,2}=  q^{12}  {\cal ABBC }
         -q^{10}  {\cal ABCB  }
         -q^7     {\cal ACAC   }
         +q^8     {\cal ACBB    }
         -q^{10}  {\cal BABC   } \]
\[       +q^8     {\cal BACB   }
         +q^8     {\cal BBAC   }
    -(q^9+q^{15}) {\cal BBBB    }
         +q^{10}  {\cal BBCA   }
         +q^{10}  {\cal BCAB   } \]
\[       -q^8     {\cal BCBA   }
         +q^{10}  {\cal CABB   }
         -q^5     {\cal CACA   }
         -q^8     {\cal CBAB   }
         +q^6     {\cal CBBA   },     \]

\[ {\cal I}_{3,3}=  q^8      {\cal AACC }
     - q^{10} [2]            {\cal ABBC }
     +      q^6              {\cal ACAC }
     - q^{10} [2]            {\cal BACB }
     + q^{12} [2]^2            {\cal BBBB } \]
\[   - q^8 [2]               {\cal BCAB }
     +      q^6              {\cal CACA }
     - q^8 [2]               {\cal CBBA }
     +      q^4              {\cal CCAA }.  \]
We have ${\cal I}_{3,1} = q^{10} {\cal I}_1^2$.
Since $F_I$ is a free algebra, the universal invariants
${\cal I}_{3,1}$,
${\cal I}_{3,2}$ and
${\cal I}_{3,3}$ are independent.

We do not try to describe the algebraic structure
of the algebra of universal invariants.
Here we consider this subject
for invariants of a certain realisation of a 2-form.

In the classical case
the algebra of invariants is generated
by the dis\-cri\-mi\-nant.
However, if $q^{2d} \neq 1$ we have two basic invariants
for $f=q(01)(02)$:

\begin{Sa} Let $q^{2d} \neq 1$ for all $d\in \NN$ and
let $f$ be the quadratic form $ q(01)(02) $.
The algebra of invariants is generated by
$I_{1,f}$ and  $I_{2,f}$.
They are related by the equation
\begin{eqnarray}
 q^2 [2] I_{2,f}^2 + (q^2-1)^2  I_{1,f}^3 = 0.
\end{eqnarray}
\end{Sa}

{\it Proof.} Let $I$ be an invariant of $f$ of degree $k$. Then
we have $I\in H_{\{1,2\}}^{Inv}$. By Theorem 3.6
we have $I=(12)^k$.
Since the symbol $(12)$ is not a linear combination
of $A,B,C$ (cf. Example 4.3), we have $k\geq 2$.
Let $k=2d+3e$, $d,e\in \NN$. By (24) and (25) we have
$I=c\  I_{1,f}^d I_{2,f}^e$
with a certain constant $c\in \CC$.
The relation (26) follows from (24) and (25).
$\bullet$

\begin{flushleft}
6.3 \it  Two and three quadratic forms
\end{flushleft}

In Example 5.4 we obtained for the symbol
$\frac{1}{q} (12)^2$ the
universal common invariant of two
forms ${\cal F}, {\cal F'}$
\[ {\cal I}_{1,{\cal F,F'}} =
            \frac{1}{q^3}  {\cal C A}'
       - (1+\frac{1}{q^2}) {\cal B B}'
       +    \frac{1}{q}    {\cal A C}'. \]

\begin{Ex} \rm
The realisations $I_{1,f,f'}$ and $I_{1,f,f'}$
are in general not equal.
For the realisation $f   = q (01)(03)$,
                    $f'  = q (02)(04)$ we obtain
\[ I_{1,f,f'} =\frac{1}{[2]}( q^5 (12)(34)-\frac{1}{q} (14)(23)  \]
and
\[ I_{1,f',f} =\frac{1}{[2]}( q^3 (12)(34)-        q^3 (14)(23)).\]
\end{Ex}

%-------------------------------------------
%\noindent
%{\it Example} 1: For the realisation $f  := q (01)(03)$,
%                   $f' := q (02)(04)$ we obtain
%\[ I_{1,f,f'}:= [2] ( q^3 Z_1 -       Z_3  )   \]
%and
%\[ I_{1,f',f}:= [2] ( q   Z_1 -   q^4 Z_3  ).   \]
%For $f  := q (01)(02)$,
%                   $f' := q (03)(04)$ we obtain
%\[ I_{1,f,f'}:= [2] ( q   Z_1 + q   [2] Z_3 )  \]
%and
%\[ I_{1,f',f}:= [2] ( q^3 Z_1 + q^3 [2] Z_3 ).  \]
%For $f  := q (01)(04)$,
%                   $f' := q (02)(03)$ we obtain
%\[ I_{1,f, f'}:= [2] (- q^2 [2] Z_1 - q^2 Z_3 )  \]
%and
%\[ I_{1,f',f }:= [2] (- q^2 [2] Z_1 - q^2 Z_3 ). \]
%
%\noindent
%{\bf Remark:} In the classical case
%5the pairs of points, corresponding to $f, f'$
%are harmonic, if the common invariant vanishes.
%--------------------------------------------------

From the symbol $(12)(01)(02)$
we obtain the common covariant of two forms
${\cal F}$, ${\cal F}'$
\[ \Delta_{ {\cal F}, {\cal F}' } =
x^2 {\cal K} + xy {\cal L} + y^2  {\cal M}  \]
with
\[ {\cal K}  = - q^7 {\cal B A}' + q^9 {\cal A B}', \]
\[ {\cal L}  = - q^7 {\cal C A}' - q^6 {\cal B B}'
            + q^{10} {\cal B B}' + q^7 {\cal A C}', \]
\[ {\cal M}  = - q^7 {\cal C B}' + q^9 {\cal B C}'. \]
In the quantised case
we have in general not the classical relation
$\Delta_{{\cal F},{\cal F}'} = -\Delta_{{\cal F}',{\cal F}}$.
%The discriminant of
%$ \Delta_{ {\cal F}, {\cal F}' }$ is the resultant of
%${\cal F}$ and ${\cal F}'$.
%The classical resultant vanishes
%if ${\cal F}$ and ${\cal F}'$ have a common zero.

\begin{Ex} \rm
{\it Example 6.5:} For  $f = f' = q(01)(02)$ we obtain
\[ \Delta_{f,f}
 = (q^\frac{23}{2} - q^\frac{15}{2}) (01)(02)(12). \]
$\Delta_{f,f}$ has the discriminant
\[ d_{f,f}=  -   q^{15} (q^2-1)^2 [2]^3 (12)^4.  \]
This discriminant
%$d_{f,f}$
vanishes only for $q^2=1$.
\end{Ex}

%For $f := q(01)(02)$
%and $f':=q(03)(04)$ we obtain
%\[ d_{f,f'} = -[2]^2 ( q^{15} (q^2-1)^2     Z_1^2    +
%                 q^{15} [2]   (1-4 q^2+q^4) Z_1 Z_3  -
%                 q^{17} [2]^2               Z_3^2   ). \]

%\noindent
%{\it Example:} $f:=q(03)(04)$, $f':=q(01)(02)$
%
%$d_{1234}=d_{3412}$
%
%\noindent
%{\it Example:} $f:=q(01)(03)$, $f':=q(02)(04)$
%
%\[  d_{1324}= - q^{12} ( q (q^4-1)^2                Z_1^2
%                     (3 + 5 q^2 + 5 q^8 + 3 q^{10}) Z_1 Z_3 +
%                       q (q^4-1)^2                  Z_3^2  ) \]
%
%\noindent
%{\it Example:} $f:=q(02)(04)$, $f':=q(01)(03)$
%
%$d_{2413}=d_{1324}$
%
%\noindent
%{\it Example:} $f:=q(01)(04)$, $f':=q(02)(03)$
%
%\[ d_{1423}:=  +  q^{13} (1 + q^2)^4               Z_1^2   -
%                  q^{12} (1 + q^2)^3 (1-4 q^2+q^4) Z_1 Z_3 -
%                  q^{13} (-1 + q^4)^2              Z_3^2    \]
%
%\noindent
%{\it Example:} $f:=q(02)(03)$, $f':=q(01)(04)$
%
%\[ d_{2314}:=  q^9 (1 + q^2)^4 (1 - q^2 + q^4)^2  Z_1^2   +
%               q^8 (1 + q^2)^3
%               (1-2 q^4 + 4 q^6 -2 q^8 + q^{12})  Z_1 Z_3 -
%               q^{13} (-1 + q^4)^2                Z_3^2   \]
%

\begin{Rem} \rm
In the classical situation
the zeros of $\Delta_{f,f'}=-\Delta_{f',f}$ are a point pair,
which is harmonic to the zeros of $f$ and to the zeros of $f'$.
Furthermore the zeros of $\Delta_{f,f'}$
are the fixpoints of the projective involution,
which is determined by the zeros of $f$ and $f'$.
We have two real zeros
if the intervalls $[v_1,v_3]$, $[v_2,v_4]$ are disjoint
or if there is an inclusion between both intervalls.
Otherwise we have two complex conjugated zeros.
This corresponds to a construction
of complex elements due to von Staudt.
In the classical situation the discriminant $d_{f,f'}$
of $\Delta_{f,f'}$ is the resultant of $f$ and $f'$,
i.e. $d_{f,f'}$ vanishes in the case of a common zero
of $f$ and $f'$.
\end{Rem}

Finally we consider an invariant of three quadratic forms
\[    f^{(i)}= x^2 A^{(i)} + [2] xy B^{(i)} +  y^2 C^{(i)},  \]
$i=1,2,3$.
From the symbol $-\frac{1}{q^4} (12)(13)(23)$ we obtain
\[{\cal I}_{f',f'',f'''} =
     q^{-\frac{  5}{2}}{\cal A'} {\cal B}'' {\cal C}'''
 -   q^{-\frac{  9}{2}}{\cal A'} {\cal C}'' {\cal B}'''
 -   q^{-\frac{  9}{2}}{\cal B'} {\cal A}'' {\cal C}'''
 +   q^{-\frac{  9}{2}}{\cal B'} {\cal C}'' {\cal A}'''
 +   q^{-\frac{  9}{2}}{\cal C'} {\cal A}'' {\cal B}'''
 -   q^{-\frac{ 13}{2}}{\cal C'} {\cal B}'' {\cal A}''' \]
\[-(q^{-\frac{7}{2}}-q^{-\frac{ 15}{2}})
                       {\cal B'} {\cal B}'' {\cal B}'''.\]
${\cal I}_{f' , f'', f'''}$ is a quantum analog of the classical
{\it combinant}
\[ \left|  \begin{array}{ccc}
               A' & A'' & A''' \\
               B' & B'' & B''' \\
               C' & C'' & C'''
\end{array}  \right|.  \]
In the classical situation
the three point pairs corresponding to $f', f'', f'''$
are in involution if and only if their combinant vanishes.

\begin{flushleft}
6.4 \it A point figuration on the real line
\end{flushleft}

%As an application of the common invariants
%$I_{1,f,f'}$, $I_{f,f',f''}$
%of two quadratic forms $f,f'$
%in order to discuss the issue
%of the harmonicity of four points.

In this section we apply the invariants and covariants
of the preceding section
to a point figuration on the real line.
In the classical case we have $ I_{1,f,f'} = I_{1,f',f}$
and the zeros of the two quadratic forms
are two point pairs, which are harmonic if and only if
$ I_{1,f,f'} = 0 $.
Furthermore we have
$ {I}_{f,f',f''} =  {I}_{f'',f,f'}  =  {I}_{f',f'',f}=
 -{I}_{f,f'',f'} = -{I}_{f',f ,f''} = -{I}_{f'',f',f} $
and the three point pairs determined by the three
quadratic forms are in involution if ${I}_{f,f',f''}=0$.
And the zeros of the quadratic form
$\Delta_{f,f'} = - \Delta_{f',f}$
are a point pair, which is harmonic
to the zeros of $f$ and $f'$.

In the following we discuss a quantum analog for points,
which are harmonic or in involution.
We consider three points $(x_i,y_i)$, $i=1,2,3$.
We form the six quadratic forms
\begin{eqnarray*}
f_{11} =  (01)(01), &  f_{22} =  (02)(02), & f_{33} =  (03)(03), \\
f_{12} = q(01)(02), &  f_{13} = q(01)(03), & f_{23} = q(02)(03).
\end{eqnarray*}
We can form the three quadratic forms
\begin{eqnarray*}
 f_{11|23}:=\Delta_{f_{11},f_{23}} - \Delta_{f_{23},f_{11}}, \\
 f_{22|13}:=\Delta_{f_{22},f_{13}} - \Delta_{f_{13},f_{22}}, \\
 f_{33|12}:=\Delta_{f_{33},f_{12}} - \Delta_{f_{12},f_{33}}.
\end{eqnarray*}
For $f_{11|23}$ and the two quadratic forms
 $f_{11}$, $f_{23}$ we have
the identities
\begin{eqnarray*}
    I_{ f_{11|23},f_{11} }=    I_{ f_{11},f_{11|23} } = 0,  \\
q^4 I_{ f_{11|23},f_{23} }+    I_{ f_{23},f_{11|23} } = 0.
\end{eqnarray*}
They give the classical harmonicity of the
corresponding forms in the classical limit.
Similar results hold for $f_{22|13}$ and $f_{33|12}$.
We have
\begin{eqnarray*}
    I_{ f_{22|13},f_{22} }=    I_{ f_{22},f_{22|13} } = 0,  \\
    I_{ f_{22|13},f_{13} }=    I_{ f_{13},f_{22|13} } = 0,  \\
    I_{ f_{33|12},f_{33} }=    I_{ f_{33},f_{33|12} } = 0,  \\
    I_{ f_{33|12},f_{12} }+q^4 I_{ f_{12},f_{33|12} } = 0.
\end{eqnarray*}
Next we consider the relation between
$f_{11|23}$, $f_{22|13}$ and $f_{33|12}$. We find
\[
  {I}_{f_{11|23},f_{33|12},f_{22|13}} =
 -{I}_{f_{22|13},f_{11|23},f_{33|12}} =
  \frac{1}{q^4} {I}_{f_{22|13},f_{33|12},f_{11|23}} =
- \frac{1}{q^4} {I}_{f_{33|12},f_{11|23},f_{22|13}}
\]
\[
=  q^7 [2]^2 (q^4+1)(q^6-1)  ((12)(13)(23))^2
\]
and furthermore
\[
{I}_{f_{11|23},f_{22|13},f_{33|12}} =  0,  \ \ \ \ \ \ \
{I}_{f_{33|12},f_{22|13},f_{11|23}} =  0,
\]
\[
{I}_{f_{33|12},f_{11|23},f_{22|13}} +
{I}_{f_{22|13},f_{33|12},f_{11|23}} = 0,
\]
\[
{I}_{f_{11|23},f_{33|12},f_{22|13}} +
{I}_{f_{22|13},f_{11|23},f_{33|12}} = 0.
\]
The formulas were checked
by a computer calculation.
Each of these four equations
gives  the classical relation
${I}_{f_{11|23},f_{22|13},f_{33|12}}=0$ in the classical limit.

\vspace{0.1cm}

This is an analog of the following classical proposition.
Let $A,B,C$ be three points and let
$P$, $Q$ and $R$
be the three fourth harmonic points of $A$, $B$ and $C$
with respect to the pair $\{B,C\}$, $\{A,C\}$ and $\{A,B\}$,
respectively.
Then the three point pairs $\{A,P\}$, $\{B,Q\}$ and $\{C,R\}$
are in involution.

\begin{flushleft}
6.5 \it  Cubic forms
\end{flushleft}

\noindent
We consider the universal cubic form
\[ {\cal F}
   =     x^3   {\cal A} + [3] x^2 y {\cal B} +
     [3] x y^2 {\cal C} + y^3       {\cal D}. \]
By  Lemma 3.3(ii) we have
\[  (ij)^3=
    q^{\frac{3}{2}}     x_i^3     y_j^3
   -q^{\frac{1}{2}} [3] x_i^2 y_i x_j y_j^2
   +q^{\frac{3}{2}} [3] x_i y_i^2 x_j^2 y_j
   -q^{\frac{9}{2}}     y_i^3     x_j^3.     \]
Applying this formula to $f=(0i)^3$, $i>0$, we obtain
\[ A= q^{\frac{3}{2}}         y_i^3,    \ \ \ \ \
   B=-q^{\frac{1}{2}}         x_i y_i^2,\ \ \ \ \
   C= q^{\frac{3}{2}}         x_i^2 y_i,\ \ \ \ \
   D=-q^{\frac{9}{2}}         x_i^3.               \]
The simplest invariant
${\cal I}_{1}$ arises from the symbol $\frac{1}{q^5} (ij)^3$.
The symbolic method gives
\[ {\cal I}_{1}=
    (-q^\frac{-19}{2} {\cal DA} + q^\frac{-13}{2} {\cal AD} )
+[3]( q^\frac{-13}{2} {\cal CB} - q^\frac{-11}{2} {\cal BC} ).\]

\begin{Ex} \rm
For the realisation $f=q^3 (01)(02)(03)$ we obtain
\begin{eqnarray}
I_{1,f}= (q^{-3}-q^3) (12)(13)(23).
\end{eqnarray}
We have $I_{1,f}\neq 0$ for $q^6 \neq 1$.
For the realisations  $f= q^2 (01)(01)(02)$,
                      $f= q^2 (01)(02)(02)$ and
                      $f=     (01)(01)(01)$
we obtain $I_{1,f}=0$.
Therefore ${\cal I}_1$ is an invariant of degree two
with the characterising property of the discriminant.
It is well known that for the classical cubic form
it is impossible to express the root of the discriminant
$(12)(13)(23)$ in terms of the coefficients.
\end{Ex}

\vspace{0.5cm}

There are no symbols,
which generate invariants of degree three.
By means of the Grassmann-Pl\"ucker relations
(cf. Proposition 3.2(iii))
we can reduce the symbols for invariants of degree four
to the following four symbols
\begin{eqnarray}
Z_1^3,     \ \ \ \ \  Z_1^2 Z_3,\ \ \ \ \
Z_1 Z_3^2,\ \ \ \ \   Z_3^3
\end{eqnarray}
(cf. Remark 3.1).
The symbols (29) can be expressed  by PBW basis elements
using a computer calculation.
Subsequently we can apply the symbolic method
to obtain the following universal invariants
${\cal D}_1, {\cal D}_2, {\cal D}_3, {\cal D}_4$:

\newpage
{\small
\[ {\cal D}_1:=
   q^{24}  {\cal  A  D  A  D }
 - q^{21}  {\cal  A  D  D  A }
 - q^{21}  {\cal  D  A  A  D }
 + q^{18}  {\cal  D  A  D  A}
\]
\[
+   [3]
(- q^{25}   {\cal  A  D  B  C }
 + q^{24}   {\cal  A  D  C  B }
 - q^{25}   {\cal  B  C  A  D }
 + q^{22}   {\cal  B  C  D  A }
\]
\[
+ q^{22}    {\cal  C  B  A  D }
- q^{19}    {\cal  C  B  D  A }
+ q^{20}    {\cal  D  A  B  C }
- q^{19}    {\cal  D  A  C  B} )
\]
\[
+(1+2 q^2+3 q^4+2 q^6+q^8)
 (  q^{22}  {\cal  B  C  B  C  }
  - q^{21}  {\cal  B  C  C  B  }
  - q^{21}  {\cal  C  B  B  C  }
  + q^{20}  {\cal  C  B  C  B}  ),
\]

\[
{\cal D}_2:=
 - q^{19}         {\cal  D  A  D  A }
 - q^{23}         {\cal  A  D  A  D}
\]
\[+ [2]
(- q^{27}         {\cal  A  C  C  C }
 - q^{23}         {\cal  C  A  C  C }
 - q^{25}         {\cal  C  C  A  C }
 - q^{21}         {\cal  C  C  C  A} )
\]
\[ + q^{28}        {\cal  A  C  B  D }
   + q^{24}        {\cal  A  C  D  B }
   + q^{23} [2]    {\cal  A  D  B  C }
   - q^{21}        {\cal  A  D  C  B}
\]
\[
 + q^{23} [2]     {\cal  B  C  A  D  }
 - q^{25}         {\cal  B  C  D  A  }
 + q^{26}         {\cal  B  D  A  C  }
 + q^{22}         {\cal  B  D  C  A }   \]
\[ + q^{24}         {\cal  C  A  B  D}
 + q^{20}         {\cal  C  A  D  B }
 - q^{21}         {\cal  C  B  A  D }
 + q^{23} [2]     {\cal  C  B  D  A }   \]
\[- q^{25}         {\cal  D  A  B  C }
 + q^{23} [2]     {\cal  D  A  C  B}
 + q^{22}         {\cal D  B  A  C }
 + q^{18}         {\cal D  B  C  A}     \]
\[ +[2]
(- q^{27}         {\cal  B  B  B  D }
 - q^{23}         {\cal  B  B  D  B }
 - q^{25}         {\cal  B  D  B  B }
 - q^{21}         {\cal  D  B  B  B} )
\]
\[
+ q^{26} [2]^2           {\cal B  B  C  C    }
- q^{21} (1+2q^2+q^4+q^{10}) {\cal B  C  B  C }
+ q^{21} [2] (1+q^8)      {\cal B  C  C  B  }
\]
\[
+ q^{21} [2] (1+q^8)      {\cal C  B  B  C  }
- q^{19} (1+q^6+2q^8+q^{10}) {\cal C  B  C  B}
+ q^{24} [2]^2           {\cal C  C  B  B},
\]

\[
{\cal D}_3:=
  q^{22}          {\cal A  D  A  D  }
+ q^{20}          {\cal D  A  D  A  }
\]
\[ + [2]
( q^{24}         {\cal A  C  C  C  }
 + q^{26}         {\cal C  A  C  C  }
 + q^{22}         {\cal C  C  A  C  }
 + q^{24}         {\cal C  C  C  A} )
\]
\[+ q^{28}          {\cal A  B  C  D }
  - q^{25}          {\cal A  B  D  C }
  - q^{26} [2]      {\cal A  C  B  D }
  - q^{21}          {\cal A  D  B  C }
\]
\[
- q^{25}         {\cal B  A  C  D }
+ q^{22}         {\cal B  A  D  C }
- q^{21}         {\cal B  C  A  D }
- q^{24} [2] {\cal B  D  A  C }
\]
\[
- q^{22} [2]  {\cal C  A  D  B }
- q^{25}          {\cal C  B  D  A }
+ q^{24}          {\cal C  D  A  B }
- q^{21}          {\cal C  D  B  A }
\]
\[
- q^{25}          {\cal D  A  C  B   }
- q^{20} [2]  {\cal D  B  C  A  }
- q^{21}          {\cal D  C  A  B }
+ q^{18}          {\cal D  C  B  A}
\]
\[ + [2]
(+ q^{24}         {\cal B  B  B  D  }
 + q^{26}         {\cal B  B  D  B  }
 + q^{22}         {\cal B  D  B  B }
 + q^{24}         {\cal D  B  B  B} )
\]
\[
- q^{22} [2] (1+q^8)       {\cal B  B  C  C }
+ q^{20} (1+q^6+2 q^8+q^{10}) {\cal B  C  B  C}
- q^{25} [2]^2            {\cal B  C  C  B  }
\]
\[
- q^{25} [2]^2            {\cal C  B  B  C }
+ q^{20} (1+2q^2+q^4+q^{10})  {\cal C  B  C  B}
- q^{20} [2](1+q^8)       {\cal C  C  B  B},
\]}

\[
{\cal D}_4:=q^{24}  {\cal A  A  D  D  }
         - q^{21}   {\cal A  D  A  D  }
         - q^{21}   {\cal D  A  D  A }
         + q^{18}   {\cal D  D  A  A} \]
\[    [3]
        (- q^{25}   {\cal A  B  C  D  }
         + q^{24}   {\cal A  C  B  D  }
         - q^{25}   {\cal B  A  D  C }
         + q^{22}   {\cal B  D  A  C}
\]
\[
         + q^{22}   {\cal C  A  D  B }
         - q^{19}   {\cal C  D  A  B }
         + q^{20}   {\cal D  B  C  A }
         - q^{19}   {\cal D  C  B  A}  )
\]
\[
+(1+2 q^2+3 q^4+2 q^6+q^8)
(          q^{22}   {\cal B  B  C  C }
 - q^{21}   {\cal B  C  B  C }
 - q^{21}   {\cal C  B  C  B }
 + q^{20}   {\cal C  C  B  B}  ).
\]

\begin{Rem} \rm
In the classical situation the discriminant
\[ d= -a^2 d^2+ 3 b^2 c^2- 4b^3 d - 4a c^3+ 6abcd \]
of the cubic form $ax^3+3bx^2y+3cxy^2+dy^3$
generates the algebra of invariants
(cf. \cite{Cl}).
If
%$q^{2d}\neq 1$ for all $d\in \NN$,
$q$ is not a root of unity
we get a four-dimensional space of invariants of degree four,
generated by
${\cal D}_1,\ {\cal  D}_2,\ {\cal  D}_3$ and ${\cal  D}_4$.
One can show that each  element of this space vanishes for
                     $f= q^2 (01)(01)(02)$,
                     $f= q^2 (01)(02)(02)$ and
                     $f=     (01)(01)(01)$.
Therefore we have a four-dimensional space of invariants
with the characterising property of the discriminant.
We have
${\cal  D}_1 = q^{37} {\cal  I}_1^2$.
In the classical case
(i.e. for $q=1$ and commuting
${\cal A},{\cal B},{\cal C},{\cal D}$)
we have ${\cal D}_2 = {\cal D}_3 = 2 d$
and ${\cal D}_1 = {\cal D}_4 = 0$.
\end{Rem}

Now we give a complete description of all invariants
for the realisation $f=q^3 (01)(02)(03)$ of $\cal F$.
The invariant ${\cal I}_1$ plays in a certain sense
the role of the root of the discriminant.

\begin{Sa}
Let $q^{2d}\neq 1$ for all $d\in \NN$
and let $f=q^3 (01)(02)(03)$.
Then the algebra of invariants is generated by $I_{1,f}$.
\end{Sa}

{\it Proof.} Let $I$ be an invariant of degree $k$.
Then we have $I\in H_{\{1,2,3\}}^{Inv}$.
By Theorem 3.6, $I$
is a linear combination of products
$(12)^{k_1} (13)^{k_2} (23)^{k_3}$
with $k_1 + k_2 + k_3= \frac{3}{2} k$.
Since $I$ is homogeneous of degree $k$
with respect to the indices $1,2,3$, we have
$ k_1+k_2 = k_1+k_3 = k_2+k_3  =  k$.
It follows that $k_1=k_2=k_3=\frac{k}{2}$.
By  Proposition 3.2 (i)
$I$ is a linear combination of products
\[ (12)^\frac{k}{2} (13)^\frac{k}{2} (23)^\frac{k}{2}
  = c ((12)(13)(23))^\frac{k}{2}
  = c'\ \ I_{1,f}^\frac{k}{2}     \]
with certain nonvanishing constants $c,c' \in \CC$,
which proves the assertion.$\bullet$

\begin{flushleft}
6.6 \it The quadratic covariant $\Delta$
\end{flushleft}

From the symbol $(12)^2 (01)(02)$
we derive  the quadratic covariant
\[ \Delta = x^2 {\cal K} +xy {\cal L} +y^2 {\cal M} \]
for the universal cubic form with
\[ {\cal K}= q^9 {\cal CA} - q^{12}[2]  {\cal BB} + q^{13} {\cal AC}, \]
\[ {\cal L}= q^9 {\cal DA} - (-q^8 +q^{12}+ q^{14}) {\cal CB} -
          ( q^9 +q^{11}- q^{15}) {\cal BC} + q^{10} {\cal AD}, \]
\[ {\cal M}= q^9 {\cal DB} - q^{12}[2]  {\cal CC}
           + q^{13} {\cal BD}.       \]
For example for $f=q^3 (01)(02)(03)$ we have
\[ \Delta_{f} = -q^{\frac{31}{2}}[2]     (01)(02)(13)(23)
                +q^{\frac{21}{2}} (1+q^8) (01)(03)(12)(23)
                -q^{\frac{23}{2}}[2]     (02)(03)(12)(13).   \]

The quadratic form $\Delta$ has the discriminant
${\cal I}_\Delta := [2]{\cal M K} - q^2 {\cal L L}
                                  + q^2 [2] {\cal K M}$ with

\[ {\cal I}_\Delta =
        - q^{20} ( q^2   {\cal ADAD}
                   + q   {\cal ADDA}
                   +     {\cal DADA}
                   + q   {\cal DAAD} )
\]
\[
 -    q^{21} [2]^2
     (q^6 {\cal ACCC} + q^2  {\cal CACC}
    + q^4 {\cal CCAC} +      {\cal CCCA} )
   + q^{24} [2]  (q^4  {\cal ACBD} + {\cal ACDB})
\]
\[
   + (q^{20} + q^{22} - q^{26})  (q {\cal ADBC}
                               +    {\cal DABC}
                               +  q {\cal BCAD}
                               +    {\cal BCDA} )
\]
\[
   +  q^{20} [2]          ( q^6   {\cal BDAC}
                          + q^2   {\cal BDCA}
                          + q^4   {\cal CABD}
                          +       {\cal CADB} )
\]
\[
   +  (-q^{19}  + q^{23}  + q^{25})
      (q {\cal CBAD} + {\cal CBDA} + q {\cal ADCB} + {\cal DACB})
   +   q^{18} [2] (q^4  {\cal DBAC} + {\cal DBCA})
\]
\[
                -   q^{21} [2]^2 ( q^6 {\cal BBBD}
                                +  q^2 {\cal BBDB}
                                +  q^4 {\cal BDBB}
                                +      {\cal DBBB} )
 +q^{24} [2]^3    (q^2 {\cal BBCC}+{\cal CCBB})
\]
\[
+(-q^{20}-2q^{22}-q^{24}+2q^{26}+2q^{28}-q^{32}){\cal BCBC}
\]
\[
+(-q^{18}+2q^{22}+2q^{24}-q^{26}-2q^{28}-q^{30}){\cal CBCB}
\]
\[
+(q^{19}+q^{21}-q^{23}-3q^{25}-q^{27}+q^{29}+q^{31}){\cal CBBC}.
\]

\begin{Rem} \rm
We have
${\cal I}_\Delta = {\cal D}_1 + [2] {\cal D}_2$.
Therefore we can distinguish the element ${\cal I}_\Delta$
in our four-dimensional space of discriminants.
In the classical case
%(i.e. for $q=1$ and commuting
%${\cal A},{\cal B},{\cal C},{\cal D}$)
we have ${\cal I}_\Delta = 2 {\cal D}_2$.
\end{Rem}

\begin{Rem} \rm
Let $f=q^3 (01)(02)(03)$ be a cubic form "with three
real zeros". We form the realisation
$\Delta_f$  of  $\Delta$.
We can consider $\Delta_f$ as an example
of a quantised real quadratic form
with two complex conjugated zeros.
The form $\Delta_f$ is related to a projective construction
of complex elements
due to F. Klein (cf. \cite{Kl}).
In the classical case
the two complex conjugated zeros
are those two unique complex points, from which one sees the
intervals, determined by the three real zeros of $f$ under an
angle of $\frac{\pi}{3}$.
%We remark, that our algebra of braided points is not compatible
%with a complex involution. (If we set $x_1:=z$, $x_2:=z^*$, we
%get an obstruction for the equation $x_2 x_1 = q^2 x_1 x_2$.)
\end{Rem}

\begin{flushleft}
6.7 \it Quartic forms
\end{flushleft}

We consider the universal quartic form
\[ {\cal F} =      x^4            {\cal A}
             +[4]  x^3 y          {\cal B}
     +\frac{[3] [4]}{[2]} x^2 y^2
                                  {\cal C} +
              [4]  x   y^3        {\cal D} +
                   y^4            {\cal E}. \]
By Lemma 3.3 (ii) we have
\begin{eqnarray}
\begin{array}{c}
   (ij)^4  =  q^4                   x_i^4     y_j^4
  -    q^2  [4]                     x_i^3 y_i x_j   y_j^3
  +    q^4  \left[  \begin{array}{c} 4 \\ 2 \end{array}   \right]
                      x_i^2 y_i^2   x_j^2 y_j^2
  -    q^4  [4]                     x_i   y_i^3   x_j^3 y_j
  +    q^8                                y_i^4   x_j^4.
\end{array}
\end{eqnarray}
Applying this formula
to the realisations $f^{(i)} =(0i)^4$, $i\neq 0$ we have
\[ A^{(i)}:=                      q^4        y_i^4, \ \ \ \ \
   B^{(i)}:= - q^2                     x_i   y_i^3, \ \ \ \ \
   C^{(i)}:=   q^2                     x_i^2 y_i^2, \]
\[ D^{(i)}:= - q^4                     x_i^3 y_i,   \ \ \ \ \
   E^{(i)}:=                      q^8  x_i^4.       \]

We consider the two symbols
\[  I_1:= \frac{1}{q^{6}} (12)^4 \ \ \ \ \ {\rm and} \ \ \ \ \
    I_2:= \frac{1}{q^{12}} (12)^2 (13)^2 (23)^2.  \]

We can use formula (30) in order to construct
the universal invariant ${\cal I}_1$ for the symbol
$I_1$,
we obtain
\[  I_1 =q^{-6}(      x_1^4          A^{(2)}
             + x_1^3 y_1   [4] B^{(2)}
             + x_1^2 y_1^2 \frac{ [3] [4] }{ [2] }
                                     C^{(2)}+
               x_1^1 y_1^3 [4] D^{(2)}+
                     y_1^4           E^{(2)}  )  \]
\[
  = q^{-6}( q^{-4}           { A}^{(1)} { E}^{(2)}-
            q^{-2} [4] { B}^{(1)} { D}^{(2)}+
            q^{-2} \frac{[3] [4]}{[2]}
                             { C}^{(1)} { C}^{(2)}-
            q^{-4} [4] { D}^{(1)} { B}^{(2)}+
            q^{-8}           { E}^{(1)} { A}^{(2)}).   \]
Therefore we have by the symbolic method
\[  {\cal I}_1   = q^{-10}   {\cal A} {\cal E} -
            q^{-8} [4]       {\cal B} {\cal D} +
            q^{-8}    \frac{[3] [4]}{[2]}
                             {\cal C} {\cal C} -
            q^{-10}[4]       {\cal D} {\cal B} +
            q^{-14}          {\cal E} {\cal A}.   \]
This is an analog of the invariant
\[ 2 (a e- 4b d+ 3 c^2)  \]
of the classical quartic form
$a x^4+4 b x^3y+ 6 c x^2y^2+4 d xy^3 + e y^4 $.
In the classical case ${\cal I}_1$ vanishes,
if the cross ratio of the four zeros
is one of the two third roots of unity
$exp( \pm \frac{2\pi {\rm i}}{3} )$,
i.e the four points are equiharmonic.

Now we construct the invariant for the symbol
$I_2= q^{-12} (12)^2 (13)^2 (23)^2$.

\newpage
A computer computation yields
\[ q^{12} I_2 = (12)^2 (13)^2 (23)^2  \]

\[ =q^4
x_1 x_1 x_1 x_1 x_2 x_2 y_2 y_2 y_3 y_3 y_3 y_3 \]
\[-( q^5  + q^{7})
 x_1 x_1 x_1 x_1 x_2 y_2 y_2 y_2 x_3 y_3 y_3 y_3 \]
\[+ q^{10}
x_1 x_1 x_1 x_1 y_2 y_2 y_2 y_2 x_3 x_3 y_3 y_3 \]
\[- (q^5  + q^{7})
 x_1 x_1 x_1 y_1 x_2 x_2 x_2 y_2 y_3 y_3 y_3 y_3 \]
\[+(- 1/q    - q + q^5  + 2 q^{7} + q^{9})
x_1 x_1 x_1 y_1 x_2 x_2 y_2 y_2 x_3 y_3 y_3 y_3 \]
\[+( q  + 2 q^3 + q^5   - q^{9} - q^{11})
x_1 x_1 x_1 y_1 x_2 y_2 y_2 y_2 x_3 x_3 y_3 y_3 \]
\[-(+ q^{7}   + q^{9})
x_1 x_1 x_1 y_1 y_2 y_2 y_2 y_2 x_3 x_3 x_3 y_3 \]
\[+ q^{10}
x_1 x_1 y_1 y_1 x_2 x_2 x_2 x_2 y_3 y_3 y_3 y_3 \]
\[+( q  + 2 q^3  + q^5  - q^{9}  - q^{11})
x_1 x_1 y_1 y_1 x_2 x_2 x_2 y_2 x_3 y_3 y_3 y_3 \]
\[+( q^{-2}  - q^2  - 3 q^4  - 3 q^{6} - q^{8}  + q^{12})
x_1 x_1 y_1 y_1 x_2 x_2 y_2 y_2 x_3 x_3 y_3 y_3  \]
\[+(- q   - q^3 + q^{7}   +2 q^{9}+ q^{11})
 x_1 x_1 y_1 y_1 x_2 y_2 y_2 y_2 x_3 x_3 x_3 y_3 \]
\[+ q^{8}
 x_1 x_1 y_1 y_1 y_2 y_2 y_2 y_2 x_3 x_3 x_3 x_3 \]
\[-( q^{7}  + q^{9})
x_1 y_1 y_1 y_1 x_2 x_2 x_2 x_2 x_3 y_3 y_3 y_3 \]
\[+(- q  - q^3  + q^{7}  +2 q^{9} + q^{11})
x_1 y_1 y_1 y_1 x_2 x_2 x_2 y_2 x_3 x_3 y_3 y_3 \]
\[+( q^3  + 2 q^5 + q^{7}  - q^{11}  - q^{13})
 x_1 y_1 y_1 y_1 x_2 x_2 y_2 y_2 x_3 x_3 x_3 y_3 \]
\[-( q^{9}  + q^{11})
x_1 y_1 y_1 y_1 x_2 y_2 y_2 y_2 x_3 x_3 x_3 x_3 \]
\[ q^{8}
 y_1 y_1 y_1 y_1 x_2 x_2 x_2 x_2 x_3 x_3 y_3 y_3 \]
\[-( q^{9}  + q^{11})
y_1 y_1 y_1 y_1 x_2 x_2 x_2 y_2 x_3 x_3 x_3 y_3 \]
\[+ q^{14} y_1 y_1 y_1 y_1 x_2 x_2 y_2 y_2 x_3 x_3 x_3 x_3. \]
We obtain by the symbolic method the universal invariant

\newpage
\begin{eqnarray}
\begin{array}{c}
{\cal I}_2 =
q^{-12} {\cal A} {\cal C} {\cal E}+
q^{-18} {\cal A} {\cal E} {\cal C}+
q^{-18} {\cal C} {\cal A} {\cal E}     \\
+q^{-16}{\cal C} {\cal E} {\cal A}+
q^{-16} {\cal E} {\cal A} {\cal C}+
q^{-22} {\cal E} {\cal C} {\cal A}     \\
+q^{-19} (1+2q^2+q^4-
q^8-q^{10})(q^2 {\cal B} {\cal C} {\cal D}+
   {\cal C} {\cal D} {\cal B}+{\cal D} {\cal B}
{\cal C})                                  \\
+q^{-21} (-1-q^2+q^6+2q^8+q^{10})(q^2 {\cal B} {\cal D}
{\cal C}+q^2 {\cal C} {\cal B} {\cal D}+{\cal D} {\cal C}
{\cal B})                                               \\
- q^{-18}[2](q^4 {\cal B} {\cal B} {\cal E}+q^2 {\cal B}
{\cal E} {\cal B}+{\cal E} {\cal B} {\cal B})   \\
- q^{-18}[2](q^4 {\cal A} {\cal D} {\cal D}+q^2 {\cal D}
{\cal A} {\cal D}+{\cal D} {\cal D} {\cal A})   \\
+ q^{-20} (1-q^4-3 q^6 -3 q^8 - q^{10} + q^{14} )
{\cal C} {\cal C} {\cal C}.
\end{array}
\end{eqnarray}
${\cal I}_2$ is an analog of the catalecticant
\[   6 \left|  \begin{array}{ccc}
               a & b & c \\
               b & c & d \\
               c & d & e
\end{array}  \right|  \]
of the classical quartic form (cf. \cite{Sa}).
In the classical case ${\cal I}_2$ vanishes,
if the cross ratio of the four zeros takes one of the values
$-1,2,\frac{1}{2}$,
i.e. the four points are harmonic.
%One can show
%that in the quantised situation
%there is no $c\in \CC$ with
%\[   {\cal I}_2 = c [12|34][13|42][14|23].  \]

In the following we will consider forms "with real zeros"
as realisations of the universal quartic form. First we consider
the form with different zeros
\begin{eqnarray}
f=q^6 (0i)(0j)(0k)(0l),\ \ \ \ \ \ \  0<i<j<k<l.
\end{eqnarray}

\begin{Sa}
The algebra of  invariants of the form $(31)$ is commutative.
\end{Sa}

{\it Proof.}
Let $I$ be an invariant of $f$ of degree $k$. By Theorem 3.6
we can represent $I$ by the bracket symbols
$(\alpha \beta )$, $\alpha, \beta =i,j,k,l,\ \alpha < \beta$.
By Proposition 3.2(i) we can reduce $I$ to a linear combination
of monomials
\[    (ij)^{k_{ij}}  (kl)^{k_{kl}}
      (ik)^{k_{ik}}  (jl)^{k_{jl}}
      (il)^{k_{il}}  (jk)^{k_{jk}}.  \]
We have
\begin{eqnarray*}
         k_{ij}+k_{ik}+k_{il} & = & 2k  \\
         k_{ij}+k_{jk}+k_{jl} & = & 2k  \\
         k_{ik}+k_{jk}+k_{kl} & = & 2k  \\
         k_{il}+k_{jl}+k_{kl} & = & 2k.
\end{eqnarray*}
If we add two of the equations and subtract the remaining two
equations in different ways, we obtain
$k_{ij}=k_{kl}$, $k_{ik}=k_{jl}$ and $k_{il}=k_{jk}$. Therefore
we can represent $I$ by a linear combination of monomials
\[    (ij)^{k_{ij}}  (kl)^{k_{ij}}
      (ik)^{k_{ik}}  (jl)^{k_{ik}}
      (il)^{k_{il}}  (jk)^{k_{il}}.  \]
If we use again Proposition 3.2(i) we can represent $I$ as
a linear combination of monomials
\[    ((ij)(kl))^{k_{ij}+i}
      ((ik)(jl))^{k_{ik}-i}
      ((il)(jk))^{k_{il}}, \ \ \ \ \ \ \ i\in \NN_0.  \]
The assertion follows by Proposition 3.2(ii).$\bullet$

Next we express the invariants
%of ${\cal I}_1$, ${\cal I}_2$
$I_{1,f}$, $I_{2,f}$ of the quartic form (32)
in terms of the "zeros".

\vspace{0.5cm}

We insert ${A}_i$ for ${\cal A}_i$ in ${\cal I}_1$
and express the resulting invariant
in terms of the bracket symbols.
By the proof of Proposition 6.4 and Remark 3.1 we can
represent the result as a linear combination of
$Z_1^2$, $Z_1 Z_3$ and $Z_3^2$.
By a computer calculation we obtain
%\[  I_1 =\frac{1}{[3][4]} (\frac{1}{q^4}
%(1-q^2+2q^3-2q^4+2q^5-q^6+q^8)
%(Z_1^2+Z_1 Z_2)
%+\frac{1}{q}(1+q^2) Z_2^2 )     \]
%or
\[ I_{1,f} = \frac{1}{[3][4]} ([2] (Z_1^2+Z_3^2) +
            (-q^{-4}+q^{-2}+2+q^2-q^4) Z_1 Z_3).  \]
Similarly we obtain
\[ I_{2,f} = \frac{[2]^2}{[3]^3 [4]^3}
   (a Z_1^3 + b Z_1^2 Z_3 + c Z_1 Z_3^2 + d Z_3^3) \]
with
\[ a=q^{-6} [2][3]
     (1+ q^2 -q^4 - q^8 + q^{10} + q^{12}), \]
\[ b=q^{-12}
        (-2+7q^4+14q^6+7q^8-8q^{10}-18q^{12}-8q^{14}
         +7q^{16}+14q^{18}+7q^{20}
         -2q^{24} ),         \]
\[ c=q^{-12} [2]
     (-1-q^2+3q^4+6q^6+2q^8-7q^{10}-13q^{12}-7q^{14}+2q^{16}
     +6q^{18}+3q^{20}-q^{22}-q^{24}),  \]
\[ d=q^{-4}  [2]^2 [3] (1-2q^2+q^4-2q^6+q^8).  \]

\begin{Rem} \rm
Using these equations we can represent the absolut invariant
$I_{2,f}^2 I_{1,f}^{-3}\in Q_I$
in terms of the quantum cross ratio $\lambda$
of the four zeros $v_i:= q^\frac{1}{2} x_i y_i^{-1}$
of the form $f$ (cf. Remark 3.2).
\end{Rem}

At the end of this section
we discuss polynomials with double zeros.
We have seven cases.

%\[ f_{iijk},\ \ f_{ijjk},\ \ f_{ijkk},\ \
%   f_{iiij},\ \ f_{iijj},\ \ f_{ijjj},\ \ f_{iiii}.  \]
%
%\[ \begin{array}{llllll}
%f_{iijk} &  0 & -q(ij)(ik) & q(ij)(ik)     &     &     \\
%f_{ijjk} &  q^2(ij)(jk) & -q^2(ij)(jk) & 0 &     &     \\
%f_{ijkk} &  0 & -q(ik)(jk) & q(ik)(jk)     &     &     \\
%f_{iiij} &  0 & 0 & 0                      & 0   & 0   \\
%f_{iijj} &  0 & -(ij)^2 & (ij)^2           &     &     \\
%f_{ijjj} &  0 & 0 & 0                      & 0   & 0.  \\
%f_{iiii} &  0 & 0 & 0                      & 0   & 0
%\end{array}   \]

\noindent
(i)For  $f_{iijk}=q^5 (0i)(0i)(0j)(0k)$
we have
\[ I_1= \frac{[2]}{[3][4]} q^4 (ij)^2 (ik)^2 \]
and
\[ I_2= \frac{[2]^4 (q^{-4}-2q^{-2}+1-2q^2+q^4)}{
      [4]^3 [3]^2 } q^9 (ij)^3 (ik)^3. \]

\noindent
(ii) For $f_{ijjk}=q^5 (0i)(0j)(0j)(0k)$
we have
\[ I_1= \frac{[2]}{[3][4]} q^8 (ij)^2 (jk)^2 \]
and
\[ I_2= \frac{[2]^3
        (-q^{-6}-q^{-4}+q^{-2}+q^2-q^4-q^6)}{
        [4]^3 [3]^2 } q^{18} (ij)^3 (jk)^3. \]

\noindent
(iii) For $f_{ijkk}=q^5 (0i)(0j)(0k)(0k)$
we have
\[ I_1=  \frac{[2]}{[3][4]} q^4 (ik)^2 (jk)^2 \]
and
\[ I_2= \frac{[2]^4 (q^{-4}-2q^{-2}+1-2q^2+q^4)}{
      [4]^3 [3]^2 } q^9 (ik)^3 (jk)^3. \]

\noindent
(iv) For $f_{iijj}=q^4 (0i)(0i)(0j)(0j)$
we have
\[ I_1= \frac{[2]}{[3][4]} (ij)^4 \]
and
\[ I_2= \frac{[2]^4 (q^{-4}-2q^{-2}+1-2q^2+q^4)}{
        [4]^3 [3]^2}    (ij)^6.    \]

\noindent
(v)-(vii)  For $f_{iiij}=q^3 (0i)(0i)(0i)(0j)$,
               $f_{ijjj}=q^3 (0i)(0j)(0j)(0j)$ and \\
               $f_{iiii}=    (0i)(0i)(0i)(0i)$
we have $I_1=0$ and $I_2=0$.

\vspace{0.5cm}

\begin{Rem} \rm
In the classical case
$6 {\cal I}_2^2 - {\cal I}_1^3$
is the discriminant of the quartic form.
In the quantised situation
the invariants ${\cal I}_1$ and ${\cal I}_2$
commute in cases (i)-(vii).
The question arises, if it is possible to combine
${\cal I}_1$ and ${\cal I}_2$
to an invariant of degree six with the
discriminant property. Up to a nonvanishing factor,
\[ {\cal D}:=
 \frac{ [3] [4]^3 }{ [2]^5 (-q^{-4}+2q^{-2}-1+2q^2-q^4)^2}
      {\cal I}_2^2-
      {\cal I}_1^3   \]
is the only element of this kind
which vanishes on\ the noncommutative polynomials
with a double zero $f_{iijj}$, $f_{iijk}$, $f_{ijkk}$,
$f_{iiij}$, $f_{ijjj}$ and $f_{iiii}$, but $\cal D$
does not vanish on $f_{ijjk}$.
There is no invariant of degree one.
${\cal I}_1$ and ${\cal I}_2$
are the only invariants of degree two and three,
respectively.
Therefore in case there exists a sixth-degree element
with the discriminant property,
there need to be further basic invariants
of degree four or six.
\end{Rem}

%
% $S:= (q^{-6}-q^{-4}-3-q^4+q^6)
%    = [3](q^{-4}-2q^{-2}+1-2q^2+q^4)$
%      \[ I_2^2- 6 I_1^3 \ \ \ \ \ \ \
%      g_2^3-27 g_3^2   \]
%

\end{document}